\documentclass[final,hidelinks,onefignum,onetabnum]{siamart220329}

\usepackage{lipsum}
\usepackage{amsfonts}
\usepackage{graphicx}
\usepackage{epstopdf}
\usepackage{amssymb,amsfonts,amsmath,latexsym,dsfont}
\usepackage{mathtools}
\usepackage{mathrsfs}
\usepackage[algo2e,ruled,vlined]{algorithm2e}
\usepackage{setspace}
\usepackage{subfigure}
\usepackage{fancyhdr}
\usepackage{wrapfig}

\usepackage{booktabs}

\usepackage[normalem]{ulem}

\usepackage{boxedminipage}
\usepackage{pgf,pgfarrows} 
\usepackage{subfigure} 

\newcommand{\FrameboxA}[2][]{#2}
\newcommand{\Framebox}[1][]{\FrameboxA}

\newcommand{\mc}[3]{\multicolumn{#1}{#2}{#3}}

\renewcommand{\div}{\nabla\cdot\,}

\newcommand{\grad}{\nabla}

\newcommand{\im}{\textit{\i}}

\newcommand{\bfe}{{\bf e}}

\newcommand{\bfp}{{\bf p}}

\newcommand{\bfu}{{\bf u}}
\newcommand{\bfq}{{\bf q}}

\newcommand{\bfr}{{\bf r}}

\newcommand{\bfmu}{{\boldsymbol \mu}}
\newcommand{\bflambda}{{\boldsymbol \lambda}}
\newcommand{\bfrho}{{\boldsymbol \rho}}

\renewcommand{\theequation}{\arabic{section}.\arabic{equation}}

\newcommand{\norm}[1]{\ensuremath{\left\|#1\right\|}}

\definecolor{darkblue}{rgb}{0.08, 0.15, 0.48}

\newtheorem{remark}{Remark}[section]

\graphicspath{{images/}{./}}

\ifpdf
  \DeclareGraphicsExtensions{.eps,.pdf,.png,.jpg}
\else
  \DeclareGraphicsExtensions{.eps}
\fi

\headers{Matrix-free multigrid for the elastic Helmholtz equation}{R. Yovel and E. Treister}

\title{LFA-tuned matrix-free multigrid method for the elastic Helmholtz equation \thanks{Corresponding author: Eran Treister. \funding{This research was supported by The Israel Science Foundation (grant No. 1589/19). RY is supported by Kreitman High-tech scholarship.}}}

\author{Rachel Yovel\thanks{Department of Computer Science, Ben-Gurion University, Beer Sheva, Israel.
  \email{yovelr@bgu.ac.il, erant@cs.bgu.ac.il}}
\and Eran Treister$^\dag$}

\usepackage{amsopn}
\DeclareMathOperator{\diag}{diag}

\ifpdf
\hypersetup{
  pdftitle={An Example Article},
  pdfauthor={D. Doe, P. T. Frank, and J. E. Smith}
}
\fi

\begin{document}

\maketitle

\begin{abstract}
We present an efficient matrix-free geometric multigrid method for the elastic Helmholtz equation, and a suitable discretization. Many discretization methods had been considered in the literature for the Helmholtz equations, as well as many solvers and preconditioners, some of which are adapted for the elastic version of the equation. However, there is very little work considering the reciprocity of discretization and a solver. In this work, we aim to bridge this gap. By choosing an appropriate stencil for re-discretization of the equation on the coarse grid, we develop a multigrid method that can be easily implemented as matrix-free, relying on stencils rather than sparse matrices. This is crucial for efficient implementation on modern hardware. Using two-grid local Fourier analysis, we validate the compatibility of our discretization with our solver, and tune a choice of weights for the stencil for which the convergence rate of the multigrid cycle is optimal. It results in a scalable multigrid preconditioner that can tackle large real-world 3D scenarios.
\end{abstract}

\begin{keywords}
Elastic wave modeling, elastic Helmholtz equation, shifted Laplacian multigrid, finite differences, high-order discretizations, full waveform inversion, parallel computations.
\end{keywords}

\begin{MSCcodes}
65N55, 74B99, 35J05, 65N06, 65F10
\end{MSCcodes}

\section{Introduction}


The Helmholtz equation, also known as the time-harmonic wave equation, is widely used for modeling propagation of waves in either acoustic or elastic media. Its applications include acoustics \cite{wang2010acoustic, elec2001acoustic}, electromagnetic radiation modeling \cite{zubeldia2012energy}, seismic modeling \cite{wang20113d} and fluid dynamics \cite{patera1984spectral, huismann2016fast}. The most common application is full waveform inversion (FWI) \cite{pratt1999, virieux2009overview, JointEikFWI17}, an inverse problem of estimating wave velocity distribution within a given domain based on observations from the boundaries. FWI is a central tool in gas and oil exploration and human brain tomography. In some of these applications the acoustic Helmholtz equation is not sufficient, and the wave propagation must be modeled by the elastic Helmholtz equation \cite{li2022target, sun2021deep, borisov20183d, kormann2017acceleration}.


Discretizing the Helmholtz equation yields a large and indefinite linear system. Both of these properties make the system difficult to solve even by sophisticated iterative methods such as multigrid \cite{briggs2000multigrid, trottenberg2000multigrid}. Moreover, modeling waves at large wavenumbers requires very fine grids. As a rule of thumb for standard finite differences discretizations, about 10 grid points per wavelength are used \cite{bayliss1985accuracy}, leading to systems that can involve hundreds of millions of unknowns. For the elastic Helmholtz equation, the linear system is even larger, for two reasons: first, it is a system of PDEs, and second, shear waves typically have a higher wavenumber compared to pressure waves, so finer meshes are required to model them \cite{martin2006marmousi2}.


When solving the equation numerically, we must take into account both the discretization and the solver. Despite the relationship between the two, in many works they are considered separately. Generally speaking, when the discretization is more sophisticated, the solver has to perform more work. A simple example of this phenomenon is the increased fill-in that occurs in direct solvers, when using high-order discretizations done by wide (non-compact) stencils. In this work we solve the elastic Helmholtz equation by a choosing the discretization in accordance with the requirements of the solution method.


There is a variety of discretizations for the Helmholtz equations. In the early \cite{babuvska1995generalized}, a finite element discretization is optimized with respect to the discretization error that is determined by the numerical dispersion. In the context of finite difference equation, the acoustic equation has the standard second-order central-difference stencil, and also compact stencils of fourth- and sixth-order \cite{singer1998high, turkel2013compact, sutmann2007compact}. These compact stencils are efficient when applying direct solvers and parallel iterative solvers on multicore processors like graphical processing units (GPUs). In the works \cite{operto20073d, aghamiry2022accurate}, combining discretizations on rotated grids to lower the numerical dispersion is discussed. It yields 9- or 27-point compact stencils with un-lumped mass terms which end up similar to \cite{singer1998high, turkel2013compact}. In \cite{stolk2016dispersion} an optimal stencil is given by a wave ray method, and \cite{stolk2014multigrid} suggests stencil coefficients that improve the multigrid convergence, based on the phase differences of the fine and the coarse grid operators. In all these high-order discretizations, typically a smaller number of grid-points per wavelength suffices.


For the elastic version, a compact nodal discretization is suggested \cite{gosselin20143d}, as well as a staggered discretization \cite{li2016Fourth, levander1988fourth}, which is more stable in the nearly incompressible case. There is also research about lowering numerical dispersion \cite{bernth2011comparison, gu201321}. However, for the staggered case the high-order discretizations are obtained by wide non-compact stencils. For example, fourth-order finite difference first derivatives in a single-dimension are used for the gradient operators in \cite{li2016Fourth, levander1988fourth}. In fact, \cite{li2016Fourth} shows that the dispersion is minimized with slightly different coefficients than the classical fourth-order terms in \cite{levander1988fourth}, hence the order of discretization is not necessarily the most important feature. Regardless, to the best of our knowledge, the available high-order schemes for the elastic equation using staggered grids are non-compact, as the extension of schemes like \cite{singer1998high, turkel2013compact} to the elastic staggered case are not straightforward \cite{li20152d}.


Concerning solvers, there are many approaches to the solution of the acoustic problem, mainly by domain decomposition \cite{gander2013domain, heikkola2019parallel} and shifted Laplacian multigrid \cite{erlangga2006novel, umetani2009multigrid, cools2014new, Dwarka2020}, but also other methods \cite{wang2020taylor, jakobsen2020convergent, azulay2022multigrid}. Fewer solvers are available for the elastic Helmholtz equation \cite{baumann2018msss, TREISTER2024112622, rizzuti2016multigrid, bonev2022hierarchical}. In this work we focus on a shifted Laplacian multigrid based solver. In our previous work \cite{TREISTER2024112622} we have shown that the mixed formulation of the elastic equation, in addition to box-smoothing, enables efficient solution of the elastic Helmholtz equation using shifted Laplacian multigrid. Also, in \cite{TREISTER2024112622} we used a standard second-order discretization accompanied with Galerkin coarsening. This coarsening strategy is natural to implement using sparse matrix computations, which have low arithmetic intensity and require excessive memory. Coarsening by re-discretization makes the solver more suitable for GPU computation and saves memory, since the stencil is explicitly given. Unfortunately, as we show and analyze in Section \ref{sec:LFA}, the multigrid cycle presented in \cite{TREISTER2024112622} shows poor convergence when using the standard second-order stencil for re-discretization. 


In this work we develop a finite difference discretization for the elastic Helmholtz equation on a staggered grid using compact stencils (i.e., using 9-point or 27-point operators in 2D and 3D, respectively). We show that our discretization is both adequate in terms of accuracy \emph{and} is suitable for stencil-based geometric multigrid to be efficient. We use mixed formulation for the elastic Helmholtz equation, and inspired by previous works \cite{singer1998high, vstekl1998accurate}, suggest a compact discretization scheme. Using local Fourier analysis (LFA), we tune the weights of the stencil such that the multigrid solver converges most efficiently, and we demonstrate that for the same weights, the numerical dispersion is low. We show that our discretization enables the geometric multigrid solver to solve the elastic Helmholtz equation even using as few as 8 grid-points per wavelength, while keeping at least the same efficiency as the standard stencil gives for 10 grid-points per wavelength.


The paper is organized as follows: in Section~\ref{sec:background} some background about multigrid methods and discretizations for the elastic Helmholtz equation is presented. In Section~\ref{sec:method} we derive our method, including the discretization scheme and the solution method. In Section~\ref{sec:LFA} we hold two-grid LFA for the system and provide theoretical results, from which the weights of the stencil are determined. We demonstrate the performance of our method in Section \ref{sec:results}, and briefly conclude in Section~\ref{sec:conclusion}.

\section{Mathematical background} \label{sec:background}

In this section we give a brief mathematical background. In Subsection~\ref{subsuc:helm} we introduce the acoustic and elastic Helmholtz equations, and derive the mixed formulation of the elastic equation. In Subsection~\ref{subsec:multigrid} we give some general multigrid preliminaries and introduce the shifted Laplacian multigrid preconditioner, and finally in Subsection~\ref{subsec:MAC} we introduce the MAC discretization scheme and appropriate multigrid components.

\subsection{The Helmholtz equations} \label{subsuc:helm}

The Helmholtz equation models wave propagation in the frequency domain, and it is in fact the Fourier transform of the wave equation. The acoustic version of the Helmholtz equation models acoustics and electromagnetic waves, whereas the elastic version models waves in solid media, such as earth's subsurface. 

Let $p = p(\vec{x}), \, \vec{x}\in\Omega$ be the Fourier transform of the wave's pressure field. The acoustic Helmholtz equation is given by:
\begin{equation}\label{eq:acousitcHelm}
\rho(\vec{x}) \div \left(\rho^{-1}(\vec{x})\grad p\right) + \kappa^2(\vec{x})  \omega^2(1-\gamma \im) p = q(\vec{x}),
\end{equation}
where $\omega = 2\pi f$ is the angular frequency, $\kappa > 0$ is the ``slowness'' of the wave in the medium (the inverse of the wave velocity), $\gamma$ represents a physical attenuation parameter and $\rho > 0$ is the density of the medium. The right-hand side $q$ represents the sources in the system. The imaginary unit is denoted by $\im.$ To solve this equation numerically, we discretize it by a finite differences scheme in a finite domain and equip it with absorbing boundary conditions (ABC) \cite{engquist1977absorbing} or perfectly matched layers (PML) \cite{berenger1994perfectly} that mimic the propagation of a wave in an open domain. 
Other boundary conditions can also be considered.

The elastic version of the Helmholtz equation has several formulations. Here we focus on the equation in an isotropic medium. Let $\vec u = \vec u(\vec x)$ be the displacement vector. The elastic Helmholtz equation is given by each of the two following formulations: 
\begin{equation}\label{eq:elasticHelm}
\grad\lambda(\vec{x})\div\vec u + \vec\div\mu(\vec{x})\left(\vec\grad\vec u+\vec\grad\vec u^T\right) +\rho(\vec{x}) \omega^2(1-\gamma \im) \vec u = \vec q_{s}(\vec{x})
\end{equation}
or equivalently\footnote{This equivalence holds for homogeneous media. In the heterogeneous case, one may use the latter as a preconditioner for the former formulation.}
\begin{equation}\label{eq:elasticHelmEquivalent}
\grad(\lambda(\vec{x}) + \mu(\vec{x}))\div\vec u + \vec\div\mu(\vec{x})\vec\grad\vec u +\rho(\vec{x}) \omega^2(1-\gamma \im) \vec u = \vec q_{s}(\vec{x}),
\end{equation}
where $\mu$ and $\lambda$ are the Lam{\'e} parameters, that define the stress-strain relationship in the media. As before, $\rho$ is the density, $\omega$ is the frequency and $\gamma$ is the attenuation. These parameters determine the pressure and shear wave velocities by $V_p = \sqrt{(\lambda+2\mu)/\rho}$, and $V_s = \sqrt{\mu/\rho}$, respectively \cite{li2016Fourth}. The term $\vec\grad\vec u+\vec\grad\vec u^T$ is the symmetric strain tensor (factored by two). The term $\vec\div \mu\vec\grad\vec u$ is the weighted diffusion operator applied on each of the components of the vector $\vec u$ separately. The Poisson ratio, defined by $\sigma = \lambda/2(\lambda+\mu)$, gives a good measure for deformation properties of the material, where most of the materials have $0<\sigma<0.5$. The nearly incompressible case, where $\sigma\to 0.5$ or equivalently $\lambda \gg \mu$, is the most difficult case for this equation as the grad-div term turns dominant.

The mixed formulation \cite{gaspar2008distributive,zhu2010efficient} is achieved by introducing a new pressure variable $p = -(\lambda+\mu)\div\vec{u}.$ The second formulation in \eqref{eq:elasticHelm}, together with the introduced pressure variable, gives the elastic Helmholtz equation the form:
\begin{equation} \label{eq:elasticMixed}
\begin{pmatrix}
\vec\div\mu\vec\grad + \rho \omega^2(1-\gamma \im) & -\grad \\
\div & \frac{1}{\lambda + \mu} Id
\end{pmatrix}
\begin{pmatrix}
 \vec u \\
 p
\end{pmatrix}
=
\begin{pmatrix}
\vec q_{s} \\
0
\end{pmatrix}.
\end{equation}

\subsection{Shifted Laplacian multigrid} \label{subsec:multigrid}

Multigrid methods \cite{brandt1977multi, trottenberg2000multigrid} are a family of iterative solvers for linear systems of the form
\begin{equation}\label{eq:linsys}
A_h\bfu=\bfq
\end{equation}
where $A_h$ is a discretized version of a given operator on a fine grid. The idea behind multigrid methods is taking advantage of the \emph{smoothing} property of standard iterative methods, such as damped Jacobi and Gauss-Seidel, to reduce the high-frequency error components, and adding a complementary \emph{coarse grid correction} process to take care of the low-frequency components. That is, we estimate and correct the error $\bfe$ for some iterate $\bfu^{(k)}$ by solving --- exactly --- a coarser analogue of the problem. To define a coarse problem, one must translate the errors from the fine grid to the coarse grid and back, using intergrid operators called \emph{restriction} and \emph{prolongation}, respectively. 

Let $P$ be the prolongation and $R$ be the restriction. Let $A_{2h}$ be the coarse grid operator, that approximates $A_h$ on the coarse grid. Then, the coarse grid correction is given by solving the equation
\begin{equation}
A_{2h}\bfe_{2h} = \bfr_{2h} = R(\bfq-A_h\bfu^{(k)})
\end{equation}
and then interpolating the solution back to the fine grid:
\begin{equation}
\bfe_h = P\bfe_{2h}.
\end{equation}
Algorithm \ref{alg:TwoCycle} summarizes this process. In matrix form, the two-grid operator is given by:
\begin{equation} \label{eq:2G}
TG = S^{\nu_2}K S^{\nu_1}
\end{equation}
where $S$ is the smoother's error propagation matrix, $\nu_1$ and $\nu_2$ are the number of pre- and post relaxations and $K$ is the coarse grid correction:
\begin{equation}\label{eq:CGC}
K = I-P A_{2h}^{-1} R A_h.
\end{equation}

\begin{algorithm}
\SetAlgoLined
\DontPrintSemicolon
\KwSty{Algorithm: $\bfu\leftarrow TwoGrid(A_h,\bfq,\bfu).$}\;
\begin{enumerate}\Indm
\item Apply $\nu_1$ times a pre-relaxation: $\bfu \leftarrow Relax(A_h,\bfu,\bfq)$\;
\item Compute and restrict the residual $\bfr_{2h} = R(\bfq - A_h\bfu)$.
\item Compute $\bfe_{2h}$ as the solution of the coarse-grid problem $A_{2h}\bfe_{2h}=\bfr_{2h}$.
\item Apply coarse grid correction: $\bfu \leftarrow \bfu + P\bfe_{2h}$.
\item Apply $\nu_2$ times a post-relaxation: $\bfu \leftarrow Relax(A_h,\bfu,\bfq)$.
\end{enumerate}
\caption{Two-grid cycle.}
\label{alg:TwoCycle}
\end{algorithm}

By applying Algorithm~\ref{alg:TwoCycle} recursively with one recursive call, we obtain the multigrid V-cycle, and with two recursive calls we obtain a W-cycle. However, as we use more levels, the ratio of frequency to coarsest grid size grows, and standard cycles tend to diverge. A nice variant is obtained by applying the two recursive calls in the W-cycle as preconditioners for two steps of a simple Krylov method \cite{notay2008recursive}. The additional cost includes two extra residual computations and a small number of vector operations on the coarse levels (e.g., inner products), which is not a lot considering the other costs. This variant is called K-cycle.

\begin{remark}
We note that the coarse grid operator $A_{2h}$ can be defined either by the Galerkin coarse approximation, as a matrix product $A_{2h}=R A_h P$ or by discretization coarse approximation, namely, re-discretizing $A_h$ (usually using the same stencil) on a coarser mesh. 
Galerkin coarse approximation has theoretical advantages as it is a projection operator, and works reasonably well in practice \cite{TREISTER2024112622}. However, discretization coarse approximation has computational advantages, since it is stencil-based, and can be easily implemented in a matrix-free manner. In this paper we use discretization coarse approximation, since we aim to develop a matrix-free solver.
\end{remark}

\paragraph*{Shifted Laplacian}
Standard multigrid methods are not effective for the solution of the acoustic Helmholtz equation~\eqref{eq:acousitcHelm}. The shifted Laplacian multigrid preconditioner suggested in \cite{erlangga2006novel} for the acoustic Helmholtz equation is based on the solution of an attenuated version of \eqref{eq:acousitcHelm} using multigrid. Let $H$ be a matrix defined by a discretization of the Helmholtz operator. We define 
\begin{equation}\label{eq:shift}
H_s = H - \im\alpha\omega^2M_s,
\end{equation}
to be the shifted Helmholtz operator, where $M_s$ is some mass matrix, and $\alpha>0$ is a shifting parameter. Adding a complex shift to the Helmholtz equation --- the wave equation in the frequency domain --- is equivalent to adding a parabolic term, that is reflected in attenuation, to the wave equation in the time domain. The added shift $\alpha$ is usually much larger than the physical attenuation $\gamma$. The shifting is implemented by adding $\alpha$ to the physical attenuation $\gamma$. For the shifted version, geometric multigrid methods are efficient, and one can use $H_s$ in \eqref{eq:shift} as a preconditioner for a discretized Helmholtz linear system \eqref{eq:linsys} inside a suitable Krylov method such as (flexible) GMRES \cite{saad1993flexible}.

\subsection{The MAC scheme for discretization} \label{subsec:MAC}

Similarly to \eqref{eq:acousitcHelm}, the elastic Helmholtz equation \eqref{eq:elasticHelm} is usually discretized using finite-differences on a regular grid. As in any system of equations, one can consider nodal discretization, in which all the variables are located in the nodes, or staggered discretization, in which every variable has a different location. The staggering is used to enhance numerical stability. For many saddle-point systems, nodal grid discretization is known to cause checkerboard instability \cite{trottenberg2000multigrid}.

\begin{figure}
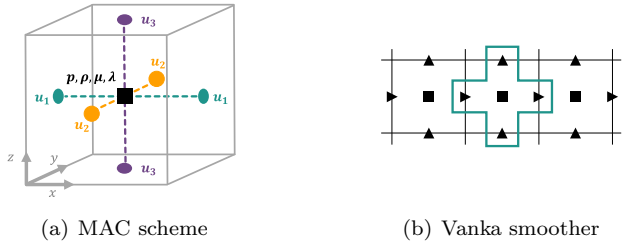

\begin{center}
	\newcommand{\image}[1]{\includegraphics[width=0.25\linewidth]{./#1}}
    \subfigure[\footnotesize MAC scheme]{\image{cell.eps}} \label{fig:stag}
    \hspace{40pt}
    \subfigure[\footnotesize Vanka smoother]{\image{vanka.eps}}\label{fig:vanka} \\
\end{center}
\caption{On the left, the MAC staggered grid discretization of a cell in 3D.
On the right, the DOFs relaxed simultaneously in one cell by the Vanka smoother in 2D.
}\label{fig:MacVanka}
\end{figure}

The MAC scheme is a common approach for a staggered grid based finite differences discretization. As depicted in Figure~\ref{fig:stag}, the displacement components are located on the faces, and the pressure is located in the center of the cells. Originally, this scheme was developed for fluid flow problems \cite{mckee2008mac, greif2022block}. It is also common for linear elasticity equation, e.g., \cite{gaspar2008distributive, he2022parameter}, and was used for the elastic Helmholtz equation in \cite{TREISTER2024112622}. To complete the discretization process, the MAC scheme should be accompanied with stencils for each component of the equation. Choosing these stencils is the core of this work, as explained later.

When using a multigrid approach to solve a problem discretized by the MAC scheme, we should define special intergrid operators as well as smoothers \cite{trottenberg2000multigrid}. Here, we elaborate only on the choices that are relevant to our work, and in two dimensions. The extension to three dimension is straight-forward.

For the pressure $p$, we take a low-order restriction and a higher-order bilinear interpolation:
\begin{equation} \label{eq:pIntergrid}
R_p=
\frac{1}{4}
\begin{bmatrix}
1 & & 1\\
 & * & \\
1 & & 1
\end{bmatrix}
\quad \text{and} \quad
P_p =
\frac{1}{16}
\left]
\begin{matrix}
1 & 3 && 3 & 1\\
3 & 9 && 9 & 3 \\
& & * & &\\
3 & 9 && 9 & 3 \\
1 & 3 && 3 & 1
\end{matrix}
\right[
\end{equation}
where [ ] represents an operator defined by a stencil, and ] [ represents the transpose of an operator defined by a stencil. The asterisk represents the center of the stencil.
Note that the pressure is cell-centered and hence the center of the stencil (where the coarse pressure is located) is not one of the sampling points. Each of the displacement components, $u_1$ and $u_2$, is located on edges on one direction and is cell-centered on the other direction, (see Fig. \ref{fig:stag}). We choose full-weighting for the first direction and, similarly to the pressure, low-order restriction and bilinear interpolation for the other direction. For the vertical component $u_2$, it reads
\begin{equation} \label{eq:u2Intergrid}
R_{u_2} =
\frac{1}{8}
\begin{bmatrix}
1 &   & 1\\
2 & * & 2\\
1 &   & 1
\end{bmatrix} 
\quad \text{and} \quad
P_{u_2} =
\frac{1}{8}
\left]
\begin{matrix}
1 & 3 &   & 3 & 1 \\
2 & 6 & * & 6 & 2 \\
1 & 3 &   & 3 & 1
\end{matrix}
\right[,
\end{equation}
and the corresponding similar operators are used for $u_1$ as well, denoted by $P_{u_1}$ and $R_{u_1}$, respectively. Finally, the restriction and prolongation operators for the whole system are defined 
\begin{equation}
P=blockDiag(P_{u_1},P_{u_2},P_p) \quad \text{and} \quad R=blockDiag(R_{u_1},R_{u_2},R_p),
\end{equation}
where $blockDiag$ forms a large block diagonal operator given the individual operators for its blocks.

As a smoother, we use the Vanka box-smoother \cite{vanka1986blockFlow}, which was originally used for fluid flow, and later adapted as a smoother for linear elasticity equation \cite{wobker2009numerical}. In this relaxation method, a whole cell is relaxed simultaneously, see Figure \ref{fig:vanka}. Namely, in each relaxation step we invert the $5\times 5$ submatrix of the fine grid operator ($7\times 7$ in 3D), that involves the DOFs of the same cell, or an approximation of this submatrix. To the version where we approximate the displacement block of the submatrix by a diagonal only, we refer as economic Vanka. This version is used for the LFA predictions and comparative results in Subsections \ref{subsec:LFAtuning} and \ref{subsec:LFApredictions}. For the results in Subsections \ref{subsec:2D} and \ref{subsec:3D} we use the full Vanka smoother, where the whole submatrix is inverted. In both full and economic Vanka, a damping parameter is used to improve the smoothing.

\section{Method} \label{sec:method}

In this section we derive our method. In \ref{subsec:disc}, we first describe a discretization method for the elastic Helmholtz equation, followed by a sketch of our multigrid cycle in \ref{subsec:cycle}.

\subsection{Discretization} \label{subsec:disc}

We observe that the main block in the mixed formulation \eqref{eq:elasticMixed}, is a block-diagonal operator with acoustic Helmholtz operators on its diagonal. Keeping in mind this observation, we discretize the acoustic blocks using an existing high-order stencil. Then, to obtain a discretization of \eqref{eq:elasticMixed}, we seek for appropriate gradient and divergence discretizations.

Our inspiration is the following stencil for the acoustic Helmholtz equation, suggested in \cite{singer1998high} as a compact fourth-order discretization for \eqref{eq:acousitcHelm} in constant coefficients:
\begin{equation} \label{eq:HO}
H^{HO} = 
\frac{1}{h^2}
\begin{bmatrix}
-1/6 & -2/3  & -1/6\\
-2/3  & 10/3 & -2/3  \\
-1/6  & -2/3  & -1/6
\end{bmatrix}
- \kappa^2 \omega^2 (1-\gamma \im)
\begin{bmatrix}
 &1/12 & \\
1/12  &  2/3 & 1/12 \\
 & 1/12  &
\end{bmatrix}.
\end{equation}
This stencil was further validated in \cite{umetani2009multigrid}, in the framework of shifted Laplacian multigrid for the acoustic Helmholtz equation. We observe that the stencil \eqref{eq:HO} can be seen as
\begin{equation} \label{eq:AcousticHbeta}
H^\beta = -\Delta_h^\beta -\kappa^2 M^\beta
\end{equation}
where $\Delta^\beta$ is a convex combination of the standard and skew Laplacian:
\begin{equation} \label{eq:LapHO}
-\Delta_h^\beta =
\beta
\frac{1}{h^2}
\begin{bmatrix}
 & -1 & \\
 -1 & 4  & -1 \\
  & -1 & 
\end{bmatrix}
+ (1-\beta)
\frac{1}{2h^2}
\begin{bmatrix}
 -1 &  & -1 \\
  & 4  & \\
 -1 &  & -1
\end{bmatrix}
\end{equation}
and $M^\beta$ is a convex combination of the identity and spread mass matrices:
\begin{equation} \label{eq:massHO}
M^\beta=
\omega^2 (1-\gamma \im) \left(
\beta\begin{bmatrix}
1 
\end{bmatrix}
+(1-\beta)\frac{1}{4}\begin{bmatrix}
& 1 & \\
1 & & 1 \\
& 1 & 
\end{bmatrix}
\right),
\end{equation}
with $\beta=2/3$ in both combinations. Our idea is to develop a parametrized discretization for the elastic Helmholtz equation in a similar manner, and tune the parameter $\beta$ such that the multigrid cycle will converge well.

\begin{figure}
  \centering
  \includegraphics[scale=0.09]{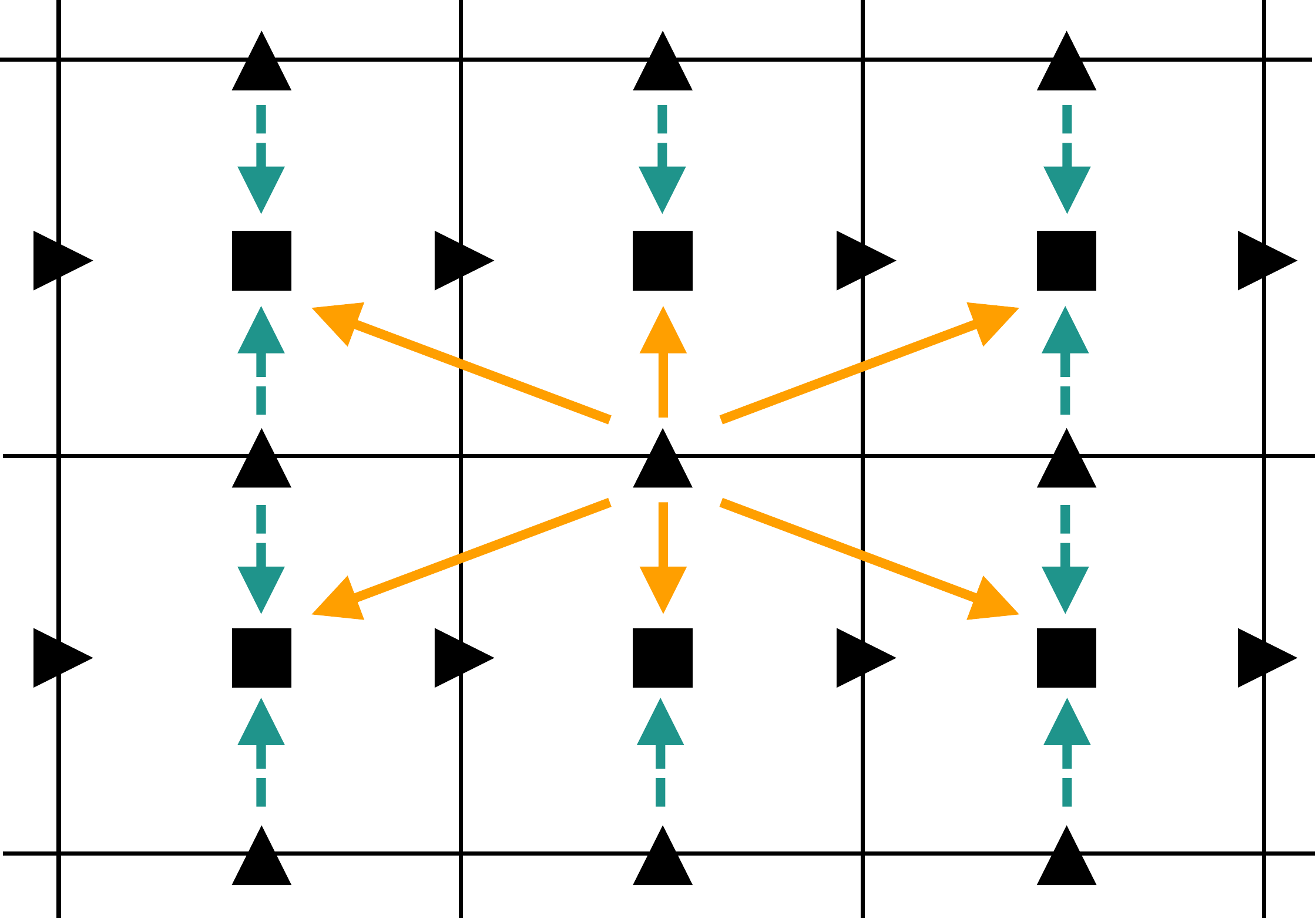}
\caption{
The sparsity patterns of the first derivative operators for the vertical component $\bfu_2$. Solid arrows denote the stencil for the spread operator $(\partial_{x_2})_{h/2}^\beta$, while dashed arrows denote the stencil for the standard operator  $(\partial_{x_2})_{h/2}$. When transposing $(\partial_{x_2})_{h/2}$ we apply the arrows on the opposite direction, and $\nabla_h^T\nabla^\beta_h = -\Delta^\beta_h$ results in a 9-point discrete Laplacian operator.}\label{fig:9ptPattern}
\end{figure}

When discretizing any of the formulations \eqref{eq:elasticHelm}, \eqref{eq:elasticHelmEquivalent} or \eqref{eq:elasticMixed} of the elastic Helmholtz equation,  discretizations of the gradient and divergence  that yield the Laplacian stencil in \eqref{eq:LapHO} by $\Delta =\div \grad$ are required\footnote{In fact, it is required for the acoustic Helmholtz equation \eqref{eq:acousitcHelm} as well, when considering non-constant coefficients.}. As depicted in Figure~\ref{fig:9ptPattern}, a compact 9-point stencil for the Laplacian can be achieved by a spread divergence based on 6-point stencils for the first derivatives, and a gradient comprised of standard 2-point stencils for the first derivatives (or vice versa).

Particularly, we discretize the horizontal first derivative as
\begin{equation} \label{eq:dxdyStandard}
(\partial_{x_1})_{h/2} = \frac{1}{h}
\begin{bmatrix}
        -1 & * & 1
\end{bmatrix}
\end{equation}
and the $\beta$-spread version as:
\begin{equation} \label{eq:dxSpread}
(\partial_{x_1})^\beta_{h/2} =
    \frac{1}{h}\left(
         \beta \begin{bmatrix}
        -1 & * & 1 \\
    \end{bmatrix}
    +
   (1-\beta)\cdot \frac{1}{4}\begin{bmatrix}
        -1 & & 1\\
        -2 & * & 2 \\
        -1 & & 1
    \end{bmatrix}
    \right),
\end{equation}
with similar (rotated) stencils for the vertical first derivatives $(\partial_{x_2})_{h/2}$ and $(\partial_{x_2})_{h/2}^\beta$. The standard and $\beta$-spread gradients are given by:
\begin{equation} \label{eq:gradSpread}
\nabla_h = 
\begin{pmatrix}
	(\partial_{x_1})_{h/2} \\[0.5em]
    (\partial_{x_2})_{h/2} 
\end{pmatrix},\quad
\nabla_h^{\beta} = 
\begin{pmatrix}
	(\partial_{x_1})^\beta_{h/2} \\[0.5em]
    (\partial_{x_2})^\beta_{h/2} 
\end{pmatrix}.
\end{equation}
The resulting Laplacian is
\begin{equation}
\grad^T_h \grad_h^\beta = -\Delta_h^\beta.
\end{equation}
For $\beta=2/3$, it is the Laplacian term in \eqref{eq:HO}.
Furthermore, the standard and $\beta-$spread divergence are given by
\begin{equation} \label{eq:divStandard}
(\nabla\cdot)_h = 
\begin{pmatrix}
	(\partial_{x_1})_{h/2} & (\partial_{x_2})_{h/2} 
\end{pmatrix}
,\quad 
(\nabla\cdot)_h^\beta = 
\begin{pmatrix}
	(\partial_{x_1})^\beta_{h/2} & (\partial_{x_2})^\beta_{h/2} 
\end{pmatrix}.
\end{equation}

Finally, we get the following discretized form of \eqref{eq:elasticMixed}:
\begin{equation} \label{eq:mixed-discretized}
\mathcal{H}^\beta 
\begin{pmatrix}
 \vec \bfu \\
 \bfp
\end{pmatrix}
=
\begin{pmatrix}
\vec\grad^T_h A_e(\bfmu) \vec\grad^\beta_h -  M^\beta A_f(\bfrho) & (\grad\cdot)_h^T \\
(\grad\cdot)_h^\beta & \diag\left(\frac{1}{\bflambda+\bfmu}\right)
\end{pmatrix}
\begin{pmatrix}
 \vec \bfu \\
 \bfp
\end{pmatrix}
=
\begin{pmatrix}
\vec \bfq_{s} \\
\mathbf{0}
\end{pmatrix}.
\end{equation}
where $A_e(\bfmu)$ is a diagonal matrix that averages the cell-centered $\bfmu$ on the edge centers, $A_f(\bfrho)$ averages cell-centered $\bfrho$ on the face centers and $\diag(\cdot)$ creates a diagonal matrix with cell-centered values on its diagonal. The shifted version $\mathcal{H}^\beta_s$ is defined by adding $\alpha$ to the physical attenuation $\gamma$.

For the acoustic Helmholtz equation, the choice $\beta=2/3$ in \eqref{eq:LapHO} and \eqref{eq:massHO} is optimal, in terms of order of the discretization, as shown in \cite{singer1998high}. In Section \ref{sec:LFA} we show, using LFA, that $\beta=2/3$ is optimal for the discretization \eqref{eq:mixed-discretized} of the elastic Helmholtz equation in mixed formulation as well, in terms of two-grid convergence.

\subsection{The multigrid cycle} \label{subsec:cycle}

Equipped with the discretization given in Section \ref{subsec:disc}, we now complete the description of our method outlining the multigrid components that we use.

As a relaxation method, we use the cell-wise Vanka smoother described in subsection \ref{subsec:MAC} and apply one pre- and one post-smoothing. For the ease of the LFA derivations in the next section, as well as for the related comparative results, we use economic Vanka in a lexicographic order. For the rest of the numerical results, we use the full Vanka smoother in red-black order to allow parallelism. This smoother was used for the elastic Helmholtz equation in \cite{TREISTER2024112622}.

As integrid operators, we use $R$ and $P$ from Section \ref{subsec:multigrid}. Note that the restriction $R$ is \emph{not} a transpose of the prolongation up to a factor. This choice is made to compensate between high-order intergrid operators, that enables a good approximation of the fine error on the coarse grid, and lower-order intergrid operators, that use smaller stencils. Finally, we define the coarse grid operator by re-discretizing the elastic Helmholtz operator using the discretization \eqref{eq:mixed-discretized}, only on a coarser grid. 

We use W-cycles on the shifted version of \eqref{eq:mixed-discretized} as a right preconditioner for the original equation inside the restarted FGMRES method. We use two, three or four levels in the multigrid hierarchy, and as we use more levels, a higher shift is required. 

\paragraph{Coarsest grid solution} The choice of the coarsest grid solver is not trivial. It can be obtained using an LU decomposition, which is our choice for 2D problems. However, this option is not practical in 3D due to high memory consumption. Another option that can be used is the domain decomposition approach in \cite{TREISTER2024112622}. This approach significantly helps memory-wise, but still includes the local LU decompositions, which are applied sequentially using forward and backward substitution and hence hinders parallelism on many-core devices like GPUs. Instead, in this work we use a parallel hybrid Kaczmarz preconditioner with GMRES, which does not require any special setup (both time and memory-wise), and can be applied in parallel. That is, we divide the domain into several sub-regions, and apply a few Kaczmarz relaxations (for each sub-region in parallel) as a preconditioner to GMRES, similarly to the approach in \cite{gordon2013robust,li20152d}. To apply it, we use the algebraic Schur complement, eliminating the $\bfp$ variable (by inverting the diagonal $\bfp-\bfp$ block), and essentially revert to the original formulation \eqref{eq:elasticHelmEquivalent}. It is done to reduce the DOFs and eliminate the need of cell-wise operations like Vanka relaxations.

\begin{remark}Note that in the case of the full Vanka relaxation, we invert the full $5\times 5$ or $7\times 7$ submatrix of each cell, for 2D or 3D, respectively. However, by the structure of \eqref{eq:mixed-discretized}, the main block that corresponds to the displacement variables is block-diagonal (containing two or three $2\times2$ sub-matrices). Thus, the submatrix that corresponds to each cell can be inverted with less operations and memory if this structure is exploited, similarly to the way the diagonal approximation is exploited in the economic Vanka variant. Utilizing the block structure, in 3D we require 25 floating numbers to store the memory for the inverted submatrix instead of 49 in the standard way. For comparison, we require 19 for the economic Vanka smoother.\end{remark}

\section{Local Fourier analysis} \label{sec:LFA}

Local Fourier analysis (LFA) is a predictive tool for the convergence of multigrid cycles \cite{brandt1977multi}. The simplest form of LFA is smoothing analysis: determining the smoothing properties of the relaxation method. Under the assumption that the coarse grid correction is ideal (a projection on the high frequencies), smoothing analysis suffices to predict the convergence rate of the two-grid cycle as a whole. 

However, sometimes the problem lies in the coarse grid correction itself, which is the case here: as shown in Subsection \ref{subsec:LFApredictions}, despite the good smoothing, the standard stencil with $\beta=1$ shows poor convergence in practice. Hence, two grid analysis gives a better prediction here.

\subsection{LFA preliminaries} \label{subsec:LFAprelim}

The definition of a smoothing factor for a system of equations (see, e.g.,~\cite{trottenberg2000multigrid}, Chapter 8) is given below.
\begin{definition} \label{def_mu_loc}
    Let $S$ be the error propagation matrix of the smoother, and let $\tilde{S}(\theta)$ be its matrix of symbols, where $\theta = \begin{pmatrix} \theta_1 & \theta_2 \end{pmatrix}^T\in\left[-\frac{\pi}{2},\frac{3\pi}{2}\right]^2$.
Then
    \[
        \mu_{loc} \coloneqq \sup_{\theta\in T^{high}} \rho(\tilde{S}(\theta))
    \]
    where $T^{high}=\left[-\frac{\pi}{2},\frac{3\pi}{2}\right]^2 \setminus \left[-\frac{\pi}{2},\frac{\pi}{2}\right]^2$.
\end{definition}
For the Vanka smoother particularly, and for overlapping smoothers generally, the calculation of the error propagation matrix requires special approaches \cite{maclachlan2011local}. For the smoothing analysis, we did computations similar to our previous work \cite{TREISTER2024112622}, only with more non-zero elements in every stencil.

The definition of the two-grid factor is given bellow:
\begin{definition} \label{def_rho_loc}
    Let $\widetilde{TG}(\theta)$ be the matrix of symbols of the two-grid operator \eqref{eq:2G}. Then
    \[
        \rho_{loc} \coloneqq \sup_{\theta\in T^{low}} \rho(\widetilde{TG}(\theta))
    \]
    where $T^{low}=\left[-\frac{\pi}{2},\frac{\pi}{2}\right]^2$. 
\end{definition}
Assuming a perfect coarse grid correction, the amplification factor of the smoother on high frequencies resembles the cycle as a whole. Namely, for an \emph{ideal} coarse grid correction, $\mu_{loc}^\nu = \rho_{loc}$, where $\nu=\nu_1+\nu_2$ is the total number of relaxation steps.

Smoothing analysis works under the assumption that the Fourier modes are eigenfunctions of the smoother. This approximately holds locally when neglecting the effect of boundary conditions. However, the Fourier modes are not eigenfunctions of the two-grid operator (not even locally) since for any low frequency $\theta\in T^{low}$, there are three high frequencies that alias to $\theta$. This gives rise to the following definition \cite{trottenberg2000multigrid}:
\begin{definition} \label{def:4harmonics}
The 4-dimensional space of harmonics for $\theta$ is
\begin{equation}
E(\theta) = \text{span} \left\{ \theta, \theta', \theta'', \theta''' \right\}
\end{equation}
where 
\[
\theta'=\begin{pmatrix} 
	\theta_1 \\ 
	\theta_2
\end{pmatrix} +
\begin{pmatrix}
\pi \\
\pi
\end{pmatrix}
,
\quad
\theta'' = \begin{pmatrix} 
	\theta_1 \\ 
	\theta_2
\end{pmatrix} +
\begin{pmatrix}
\pi \\
0
\end{pmatrix}
\quad \text{and} \quad
\theta''' = \begin{pmatrix} 
	\theta_1 \\ 
	\theta_2
\end{pmatrix} +
\begin{pmatrix}
0 \\
\pi
\end{pmatrix}
.
\]

\end{definition}

The symbol matrix of the two-grid cycle is:
\begin{equation} \label{eq:symbolTG}
\widetilde{TG}(\theta) = \tilde{S}(\theta)^{\nu_2} (I-\tilde{P}(\theta) \tilde{A}_{2h}^{-1}(\theta) \tilde{R}(\theta) \tilde{A}_h(\theta)) \tilde{S}(\theta)^{\nu_1}
\end{equation}
where the dimensions of each symbol matrix are determined by the space of harmonics. For the acoustic Helmholtz equation, a detailed description is given in \cite{cools2013local}. As mentioned there, once we calculate the symbol of the smoother as a function of $\theta$, we need to evaluate it in these four basis elements of the space of harmonics, and place it on a diagonal a $4\times 4$ matrix of symbols. In a similar manner, the symbol matrix of the fine grid operator is calculated. The restriction and prolongation symbol matrices are $1\times 4$ and $4 \times 1$ respectively, and the symbol of the coarse operator is a scalar.

For our system of equations, the dimensions of the symbol matrices must be larger in correspondence to the number of variables. In 2D, we repeat this process for the variables $u_1,u_2,p$ and get a $12\times 12$ matrix for the symbol of the two-grid operator. The fine operator's matrix of symbols, as well as the smoother's, is $12\times 12$, the matrix of symbols for the restriction and prolongation are $3\times 12$ and $12 \times 3$ respectively, and the symbol of the coarse grid is a $3\times 3$ matrix.

\subsection{Derivation of the two-grid symbol}

In this subsection we derive two-grid LFA for the discretization \eqref{eq:mixed-discretized} in 2D. As LFA is local, we can only predict convergence for the case of constant coefficients, for which \eqref{eq:mixed-discretized} takes the form:
\begin{equation} \label{eq:mixed-discretized-homogeneous}
\begin{pmatrix}
 -\bfmu \Delta^\beta_h - \bfrho M^\beta & & (\partial_{x_1})_{h/2} \\[0.15em]
& -\bfmu \Delta^\beta_h - \bfrho M^\beta &  (\partial_{x_2})_{h/2} \\[0.15em]
(\partial_{x_1})_{h/2}^\beta & (\partial_{x_2})_{h/2}^\beta & \frac{1}{\bflambda+\bfmu}I
\end{pmatrix}.
\end{equation}

In order to calculate the symbol $\widetilde{TG}(\theta)$ of the two-grid operator \eqref{eq:2G}, we calculate the symbol matrix of the smoother and the coarse grid correction. 

The Vanka smoother is an overlapping smoother, as some of the DOFs are corrected twice in each swap. The calculation of the symbol for overlapping smoothers is very lengthy, see \cite{maclachlan2011local}. In our previous work \cite{TREISTER2024112622}, we gave the analysis for the case of standard stencils (5-point stencil for the acoustic Helmholtz block and 2-point stencils for the first derivatives). The generalization for a 9-point stencil for the acoustic block and 6-point stencil for the gradient's first derivatives is straightforward, and we omit the gory details. We denote by $\widetilde{S}_h(\theta)$ the calculated scalar symbol of the Vanka smoother, and construct the smoothing matrices as following:
\begin{equation} \label{eq:smoothingSymbolMatrix}
\widetilde{S} = Diag(\widetilde{S}(\theta), \widetilde{S}(\theta'), \widetilde{S}(\theta''), \widetilde{S}(\theta''')).
\end{equation}
For the coarse grid correction, we give a more detailed analysis of the symbol.

To construct the symbol matrices of the restriction and the prolongation, we first compute the scalar symbol of each of its components. The scalar symbol of each of $R_{u_1},R_{u_2},R_p,P_{u_1},P_{u_2}$ and $P_p$ (whose stencils are given in \eqref{eq:u2Intergrid} and \eqref{eq:pIntergrid}) is calculated directly from the stencil. For instance,
\begin{align} \label{eq:Ru1Symbol}
\widetilde{R}_{u_1} (\theta)& = e^{\im \left(\frac{\theta_2}{2}-\theta_1\right)} + e^{\im \left( \frac{\theta_2}{2}+\theta_1 \right)} + e^{\im \left( \frac{-\theta_2}{2}+\theta_1 \right)} + e^{\im \left( \frac{-\theta_2}{2}-\theta_1 \right)} + 2e^{\im \frac{\theta_2}{2}} 2e^{- \im \frac{\theta_2}{2}} \\ \nonumber
& = 2\cos(\theta_1 +\theta_2 /2) + 2\cos(\theta_1 - \theta_2 / 2) + 4\cos(\theta_2 / 2).
\end{align}
Applying it on each of the four harmonics gives, for the restriction, the vector 
\begin{equation} \label{eq:Ru1Harmonics}
\widetilde{R}_{u_1} = \begin{pmatrix}\widetilde{R}_{u_1}(\theta) & \widetilde{R}_{u_1}(\theta') & \widetilde{R}_{u_1} (\theta'') & \widetilde{R}_{u_1} (\theta''')\end{pmatrix}.
\end{equation}
For the prolongation, we construct a similar vector with transposed dimensions:
\begin{equation} \label{eq:Pu1Harmonics}
\widetilde{P}_{u_1} = \begin{pmatrix}\widetilde{P}_{u_1}(\theta) & \widetilde{P}_{u_1}(\theta') & \widetilde{P}_{u_1} (\theta'') & \widetilde{P}_{u_1} (\theta''')\end{pmatrix}^T.
\end{equation}

In the same way, $\widetilde{R}_{u_2}$, $\widetilde{R}_p$, $\widetilde{P}_{u_2}$, $\widetilde{P}_p$ are calculated, and finally, the symbol matrices of the intergrid operators are given by 
\begin{equation} \label{eq:symbolRestriction}
\widetilde{R} = blockDiag(\widetilde{R}_{u_1}, \widetilde{R}_{u_2}, \widetilde{R}_{u_p})
\quad \text{and} \quad
\widetilde{P} = blockDiag(\widetilde{P}_{u_1}, \widetilde{P}_{u_2}, \widetilde{P}_{u_p}).
\end{equation}

The symbol matrix of the coarse grid operator, $\widetilde{\mathcal{H}}_H,$ is a $3\times 3$ matrix comprised of the symbols of each block in \eqref{eq:mixed-discretized-homogeneous}. The symbol matrix of the fine grid operator $\widetilde{\mathcal{H}}_h,$ is a $12\times 12$ matrix, comprised of $4\times 4$ diagonal blocks: each of the scalar symbols of the blocks in \eqref{eq:mixed-discretized-homogeneous} is evaluated on each of the four harmonics from Definition \ref{def:4harmonics}, to form the diagonal blocks.

Finally, the symbol matrix of the two-grid cycle is calculated by \eqref{eq:symbolTG}.

\subsection{Tuning the stencil by LFA}\label{subsec:LFAtuning}

In this subsection we tune the discretization suggested in \eqref{eq:mixed-discretized}, namely, determine the optimal $\beta$ in \eqref{eq:mixed-discretized-homogeneous} using two-grid LFA. We do it by numerically searching for a value that minimizes $\rho_{loc}$ from Definition \ref{def_rho_loc}. Our default parameters are density $\rho=1$ and Lam{\'e} coefficients\footnote{The Lam{\'e} parameters $\lambda$ and $\mu$ has units of pressure. However, their dimensionless ratio uniquely determines the Poisson ratio $\sigma$, so for simplicity, we omit their units throughout the numerical results.} $\lambda=500$, $\mu=1$. For this experiment we take a grid size of $h=1/1024$, frequency $\omega=\pi/5h$ that corresponds to 10 grid-point per wavelength, and shift $\alpha=0.1$. To choose the damping parameter $w$, we use smoothing analysis, since we observe that the optimal damping for the two-grid cycle has a strong dependence on the choice of the default parameters. We estimate the smoothing factor $\mu_{loc}$ from definition \ref{def_mu_loc} by sampling the Vanka smoother symbol for $\theta\in[-\pi/2,3\pi/2]$ with jumps of $0.01$ in each component, and then taking the maximum over $\theta\in T^{high}$.

\begin{figure}
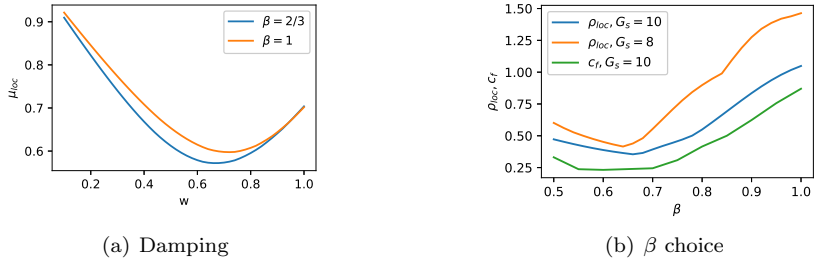

\begin{center}
	\newcommand{\image}[1]{\includegraphics[height=0.25\linewidth]{./#1}}
    \subfigure[\footnotesize Damping]{\image{muloc_w.eps}\label{fig:muloc_w}} 
\hspace{40pt}
    \subfigure[\footnotesize $\beta$ choice]{\image{rholoc_cf_beta.eps}\label{fig:rholoc_beta}} \\
\end{center}
\caption{\footnotesize Tuning the stencil. On the left, the smoothing factor as a function of the damping $w$, for grid size $h=1/1024$, frequency $\omega=\pi/5h$ that corresponds to 10 grid points per wavelength, Lam{\'e} coefficients $\lambda=500$ and $\mu=1$, density $\rho=1$ and a shift of $\alpha=0.1$. On the right, the theoretical two-grid factor $\rho_{loc}$, and the convergence factor in practice $c_f$ as a function of $\beta$. We use the above-mentioned parameters as well as an additional frequency $\omega=\pi/4h$ that corresponds to 8 grid points per wavelength, and a damping of $w=0.7$ (which is approximately optimal for the range of $\beta$ we examine). Evidently, $\beta=2/3$ is optimal for 10 grid points per wavelength, and nearly optimal for 8 grid points per wavelength.}
\label{fig:tuning}
\end{figure}

Figure \ref{fig:muloc_w} shows that the smoothing factor $\mu_{loc}$ is approximately minimized when $w=0.75$, for $\beta=1$, and when $w=0.65$, for $\beta=2/3$ (in fact, it is approximately minimized for all values of $\beta$ when $w=0.7$). The attenuation has almost no effect on the smoothing factor \cite{TREISTER2024112622}, and the smoothing factors are $\mu_{loc}^{opt}=0.55$ for $\beta=2/3$ and $\mu_{loc}^{opt}=0.59$ for $\beta=1$. Consequently, smoothing analysis is not delicate enough to distinguish between different choices of $\beta$, and two-grid analysis is necessary to tune our discretization.

Figure \ref{fig:rholoc_beta} shows the tuning of the stencil. We estimate the two-grid factor $\rho_{loc}$ with 1 pre- and one post smoothing, by sampling the two-grid symbol over $T^{low}$ in jumps of $0.01$ and then taking the maximum. With the same default parameters and damping of $w=0.7$, we consider values of $0.5<\beta<1$ with jumps of $0.02$. We observe numerically that the minimal two-grid factor $\rho_{loc} = 0.38$ is achieved when $\beta_{opt}=0.667\approx 2/3$ for the frequency that corresponds to 10 grid points per wavelength. We also observe that the optimal value of $\beta$ slightly decreases when the frequency increases. Finally, we compare these results to the convergence factor in practice, measured as in \eqref{eq:convFactor} for $0.5<\beta<1$ with jumps of $0.05$, and we see that the optimal convergence factor $c_f = 0.23$ is achieved when $\beta=0.6$, although the convergence factors for the range $0.55<\beta<0.7$ are nearly the same. To summarize, for all the frequencies that we use in practice, $\beta=2/3$ is either optimal or nearly optimal, both in the theoretical analysis and in practice. 

\begin{figure}
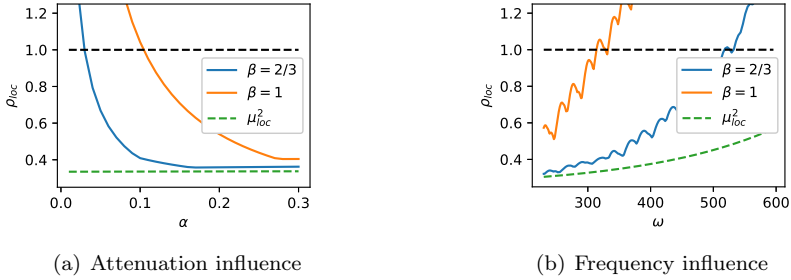

\begin{center}
	\newcommand{\image}[1]{\includegraphics[width=0.35 \linewidth]{./#1}}
    \subfigure[\footnotesize Attenuation influence]{\image{minimal_alpha_small.eps} \label{fig:min-alpha}} 
    \hspace{40pt}
    \subfigure[\footnotesize Frequency influence]{\image{maximal_omega_small.eps} \label{fig:max-omega}}  \\
\end{center}
\caption{\footnotesize On the left, the two-grid factor $\rho_{loc}$ as a function of the attenuation $\alpha$, for $\omega=\pi/5h$ which corresponds to 10 grid-points per wavelength. On the right, $\rho_{loc}$ as a function of the frequency $\omega$, for $\alpha=0.1$. For both, we use grid size $h=1/1024$, Lam{\'e} coefficients $\lambda=500$ and $\mu=1$, and density $\rho=1$.}
\label{fig:min-alpha-max-omega}
\end{figure}

Next, we investigate the influence of the shift and the number of grid points per wavelength on the convergence. Figure \ref{fig:min-alpha} shows the two-grid factor as a function of the shift $\alpha$, for the above-mentioned default parameters. We observe that for a frequency that correspond to 10 grid-points per wavelength, when $\beta=2/3$, an attenuation of $\alpha_{min}=0.03$ suffices for convergence, whereas for $\beta=1$, a larger shift of $\alpha=0.12$ is required. In practice, when the shift is smaller, the preconditioner is closer to the original system, which may improve the GMRES performance.

Figure \ref{fig:max-omega} shows $\rho_{loc}$ as a function of $\omega$, the angular frequency\footnote{The jumps in the graph occur since for some values of $\omega$, the operator has a zero eigenvalue.}. When $\beta=2/3$ and $\alpha=0.1$, the two-grid cycle still converges for about $\omega_{max}=500$, which corresponds to about 6.5 grid-points per wavelength. For comparison, using the same shift, the standard discretization requires 11 grid points per wavelength to converge. 

In the next section we compare between the discretization \eqref{eq:mixed-discretized-homogeneous} with $\beta=1$ (namely, the standard 5-point stencil), and with $\beta=2/3$.

\section{Numerical Results} \label{sec:results}

In this section we show results for 2D and 3D problems. Some of our examples are based on geophysical models, in which the length of the domain is high (about 20km) compared to its depth (about 5km). We take the right-hand side $\vec\bfq$ to be a point source located at middle of the top row (or, in 3D, the top surface) of the domain. In all the examples (except, of course, the LFA), we use an absorbing boundary layer of 20 cells, and assume that the unpreconditioned equation \eqref{eq:elasticMixed} has a small physical attenuation of $\gamma = 0.01$.

In Subsection \ref{subsec:LFApredictions} we use a two-level $V(1,1)$-cycle, with Vanka smoother in lexicographic order, to compare the actual convergence rate of the cycle with the LFA predictions. In Subsections \ref{subsec:2D} and \ref{subsec:3D}, we use the preconditioned GMRES(5) Krylov solver \cite{saad1993flexible}, preconditioned by $W(1,1)$-cycles solving the shifted version with Vanka cell-wise relaxation applied in red-black order. We seek a solution with relative residual accuracy of $10^{-6}$, starting from a zero initial guess.

We write our code in the {\tt Julia} language \cite{Julia}, and include it as a part of the {\tt jInv.jl} package \cite{jInv17}. This package enables the use of our code as a forward solver for three-dimensional elastic full waveform inversion in the frequency domain. We compute the tests on a workstation with Intel Xeon Gold 5117 \@ 2GHz X 2 (14 cores per socket) with 256 GB RAM, that runs on Centos 7 Linux distribution.

\paragraph{Demonstration of Numerical Dispersion}

\begin{figure}
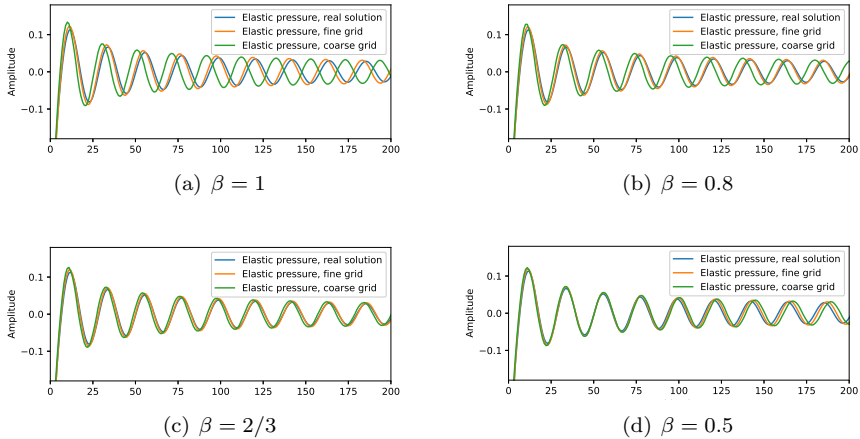

\begin{center}
	\newcommand{\image}[1]{\includegraphics[width=0.45\linewidth]{./#1}}
    \subfigure[\footnotesize $\beta=1$]{\image{beta1.eps} \label{fig:dispersion_1}} 
	\subfigure[\footnotesize $\beta = 0.8$]{\image{beta08.eps} \label{fig:dispersion_08}} 
    \subfigure[\footnotesize $\beta = 2/3$]{\image{beta066.eps} \label{fig:dispersion_066}} 
    \subfigure[\footnotesize $\beta = 0.5$]{\image{beta05.eps} \label{fig:dispersion_05}}  \\
\end{center}
\caption{\footnotesize A section of the real part of the pressure component of the direct solution, with $\lambda=\rho=16$ and $\mu=1$, $\omega=50\pi$ and grid sizes $1024\times512$, $512\times256$ and $256\times128$ for the real, fine and coarse operators respectively, for different choices of $\beta$.}
\label{fig:dispersion}
\end{figure}

\begin{figure}
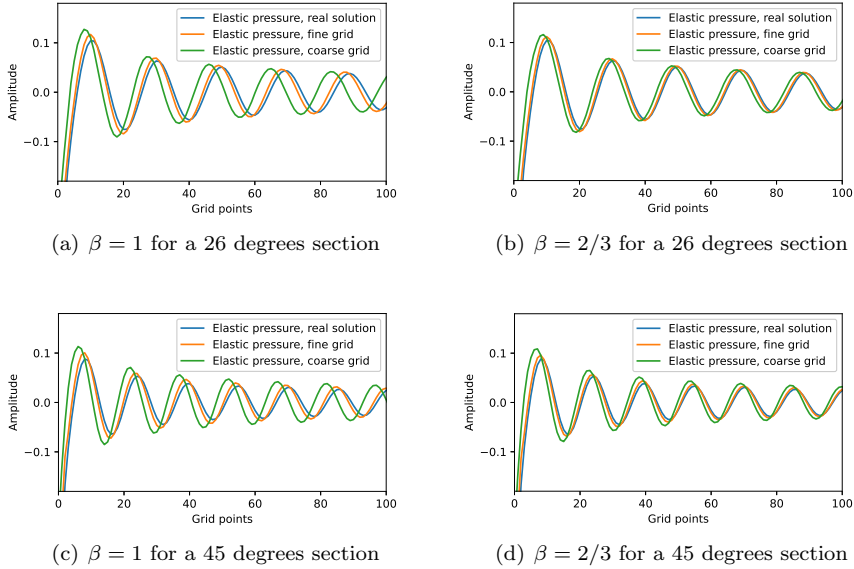

\begin{center}
	\newcommand{\image}[1]{\includegraphics[height=0.25\linewidth]{./#1}}
    \subfigure[\footnotesize $\beta=1$ for a $26$ degrees section]{\image{dispersionStandard26deg.eps} \label{fig:dispersion_standard_26}} 
	\subfigure[\footnotesize $\beta = 2/3$ for a $26$ degrees section]{\image{dispersionSpread26deg.eps} \label{fig:dispersion_spread_26}} 
    \subfigure[\footnotesize $\beta = 1$ for a $45$ degrees section]{\image{dispersionStandard45deg.eps} \label{fig:dispersion_standard_45}} 
    \subfigure[\footnotesize $\beta = 2/3$ for a $45$ degrees section]{\image{dispersionSpread45deg.eps} \label{fig:dispersion_spread_45}}
\end{center}
\caption{\footnotesize Different sections of the real part of the pressure component of the direct solution, with $\lambda=\rho=16$ and $\mu=1$, $\omega=50\pi$ and grid sizes $1024\times512$, $512\times256$ and $256\times128$ for the real, fine and coarse operators respectively.}
\label{fig:dispersionAngles}
\end{figure}

First, we demonstrate the effect of our stencil on the numerical dispersion. We note that, as mentioned in \cite{stolk2014multigrid}, the effectiveness of the coarse grid correction can be measured by the numerical dispersion between the fine and coarse grid solutions, hence this demonstration gives another justification for our choice of a stencil. For the demonstration, we solve the 2D homogeneous elastic Helmholtz equation in mixed formulation by a direct inversion of the matrix. We solve the equation on the dimensionless domain $[0,2]\times[0,1]$ with frequency $\omega=50\pi,$ density $\rho = 16$ and Lam{\'e} parameters $\lambda=16$ and $\mu=1$. We use an absorbing boundary layer of 20 grid points. Figure \ref{fig:dispersion} depicts a vertical section of the real part of the pressure component of the solution, starting from the middle of the top row, where the point source is located. We compare between the solution for 24 grid points per wavelength, which we refer to as the ``real'' solution obtained on a grid sized 1024$\times$512, compared to the solutions for 12 and 6 grid points per wavelength, referred to as the fine and coarse grid solutions respectively, obtained on grids of sizes 512$\times$256 and 256$\times$128 respectively. Figure \ref{fig:dispersion_1} demonstrates that when discretizing the equation with the standard discretization, the fine grid operator represents a significantly dispersed wave compared to the real solution, and the coarse grid solution is not only dispersed relatively to the real solution, but also fails to resemble the fine grid solution. Figure \ref{fig:dispersion_066} shows that the numerical dispersion is much lower when using our stencil. Although our choice of $\beta$ is optimal in terms of two-grid analysis, analytically determining its optimality in terms of numerical dispersion, is beyond the scope of this work. However, we show two more values of $\beta$ in Figure \ref{fig:dispersion_08} and Figure \ref{fig:dispersion_05}, and indeed, $\beta=2/3$ is optimal out of the values we examined. Moreover, we comment that the dispersion can be different when taking sections in different angles between $0$ and $45$ degrees with respect to the above mentioned vertical section. In Figure \ref{fig:dispersionAngles} we show the results for two more sections in different directions, and even though the dispersion is slightly larger for $45$ degrees than it is for a vertical section, it is still much lower for $\beta=2/3$ compared to $\beta=1$, for any of the examined sections.

\subsection{LFA predictions vs. multigrid performance}
\label{subsec:LFApredictions}

In this subsection we show the expected multigrid performance of our discretization in the 2D homogeneous case \eqref{eq:mixed-discretized-homogeneous} compared to the actual multigrid performance. We use the smoothing factor $\mu_{loc}$ from Definition \ref{def_mu_loc} to predict the convergence rate of the multigrid cycle assuming an ideal coarse grid correction. More explicitly, for a cycle with a total number of $\nu$ pre- and post smoothing steps, we use $\mu_{loc}^\nu$ as a measure for the best possible convergence rate that we can hope for in our problem. In our two-grid LFA, we use the two-grid factor $\rho_{loc}=\rho_{loc}(\nu)$, from definition \ref{def_rho_loc}, to predict the convergence of the two-grid cycle.  

Finally, we compare the results to the convergence factor in practice, defined
\begin{equation} \label{eq:convFactor}
c_f^{(k)} = \left(\frac{\norm{r_k}}{\norm{r_0}}\right)^{1/k}
\end{equation}
where $r_0$ is the residual in the error-residual equation of the two-grid operator \eqref{eq:2G} after a warm-up of 5 iterations, and $r_k$ is the residual after $k$ more iterations. We take $k$ to be the smallest number of iterations such that $r_k<10^{-9}$.

In Table \ref{tab:rho_cf_comparison} we compare the values of $\rho_{loc}$ and $c_f$ for different values of frequency $\omega$ and shift $\alpha$. We choose the rest of the parameters as our default parameters as chosen in section \ref{subsec:LFAtuning}. To avoid a significant influence of the boundaries, we take a large enough grid, with $h=1/1024$. As a reference value for the best case scenario, we use $\mu_{loc}^2$ that corresponds to a two-grid cycle with one pre- and one post-smoothing, assuming an ideal coarse grid correction.

\begin{table}
\centering
\begin{tabular}{c|ccc|ccc|ccc}
\toprule
  \mc{10}{c}{LFA two-grid factor and convergence factor in practice}\\
 \midrule
  & \mc{3}{c|}{$\omega=\frac{\pi}{5h},\,\alpha=0.15$} &  \mc{3}{c|}{$\omega=\frac{\pi}{4h},\,\alpha=0.2$} &  \mc{3}{c}{$\omega=\frac{\pi}{3.3h},\,\alpha=0.3$} \\
 discretization &  $\rho_{loc}$ & $c_f$ & $\mu_{loc}^2$ & $\rho_{loc}$ & $c_f$ & $\mu_{loc}^2$ & $\rho_{loc}$ & $c_f$ & $\mu_{loc}^2$ \\
\midrule
$\beta=1$   & 0.73 & 0.6 & 0.35 & 0.81 & 0.74 & 0.38 & 0.75 & 0.7 & 0.44 \\
$\beta=2/3$ & 0.37 & 0.24 & 0.31 & 0.4 & 0.27 & 0.34 & 0.47 &  0.35 & 0.39 \\ 
\bottomrule
 \end{tabular}
\caption{The LFA two-grid factor $\rho_{loc}$ and the convergence factor in practice $c_f$ for a grid of $h=1/1024$ with $\lambda=500$, $\mu=1$ and $\rho=1$. For the standard discretization we take damping of $w=0.75$ and for the spread discretization $w=0.65$. As a reference value, $\mu_{loc}^2$ resembles a two-grid with 1 pre- and 1 post-smoothing, assuming an ideal coarse grid correction.}
\label{tab:rho_cf_comparison}
\end{table}

\subsection{2D experiments} \label{subsec:2D}

In this subsection we demonstrate our approach, providing numerical results for several 2D models. In all experiments, we solve \eqref{eq:elasticMixed} with physical attenuation of $\gamma = 0.01$ using FGMRES and use a W(1,1)-cycle for the shifted version as a preconditioner, with varying choices of shift $\alpha$, depending on the setup (e.g., number of levels). We compare our discretization --- namely, \eqref{eq:mixed-discretized} with $\beta=2/3$ --- and the standard discretization with $\beta=1$, in terms of iterations count to reach the convergence criteria. We denote the number of grid-points per shear wavelength by $G_s$. In all the experiments we use $G_s=10$ for the standard discretization, and compare it with $G_s=10$ and $G_s=8$ using our discretization with $\beta=2/3$. For $G_s=10$, we choose damping parameters of $w=0.55,0.35,0.25$ for 2-, 3-, and 4-level respectively. For $G_s=8$, we take $w=0.5,0.35,0.25$. These values were chosen based on the performance for our red-black relaxation, and give a reasonable convergence for all the heterogeneous cases we examine. Our choice is also inspired by \cite{calandra2013improved}, which suggests for the acoustic Helmholtz equation \eqref{eq:acousitcHelm} the use of a smaller damping parameter for coarser grids. Note that the same grid with less grid points per wavelength behaves like a partially coarsened version of the original grid.

In the first experiment we use homogeneous media model: we solve \eqref{eq:elasticMixed} with constant coefficients $\lambda = 20$, $\mu=\rho=1$. We examine different grid sizes. In the case of two-grid cycles, we observe that our method converges significantly better with $\beta=2/3$ than with $\beta=1$, taking 10 grid points per wavelength, even when using a lower shift. Moreover, the performance of our stencil with only 8 grid points per wavelength is comparable to the performance of the standard stencil with 10 grid points per wavelength. The results are summarized in Table \ref{tab:2Dhomogeneous}.

\begin{table} 
\centering
\begin{tabular}{c|ccc|ccc}
\toprule
  \mc{7}{c}{Iteration count for homogeneous media}\\
\midrule
 &  \mc{3}{c|}{$\beta=1$, $G_s=10$} & \mc{3}{c}{$\beta=2/3$, $G_s=10 \, (8)$} \\
\midrule
 & \small 2-level & \small 3-level & \small 4-level & \small 2-level  & \small 3-level & \small 4-level \\
Grid size & {\small $\alpha = 0.1$} &  {\small $\alpha=0.4$}  & \small $\alpha=0.5$   & \small $\alpha=0.1$   & \small $\alpha=0.3 \,(0.5)$ & \small $\alpha=0.4 \,(0.5)$\\
\midrule
$512\times 128$ &  39 &  131 & 160 & 22 (27) & 86 (144)  & 143 (171) \\
$1024\times 256$ & 56  & 199 & 243 & 32 (39) & 135 (231) & 222 (269) \\
$2048\times 512$ & 79  & 256 & 332 & 40 (50) & 144 (294) & 275 (294) \\
\bottomrule
 \end{tabular}
\caption{Number of preconditioning cycles needed for convergence with shifted Laplacian multigrid for the 2D elastic Helmholtz equation with constant coefficients: $\mu=\rho=1$, and $\lambda=20$. The frequency corresponds the number of grid-points per shear wavelength, $G_s$. The value $\lambda=20$ corresponds to a Poisson’s ratio $\sigma = \lambda/2(\lambda+\mu)=0.476$.}
\label{tab:2Dhomogeneous}
\end{table}

In the second experiment, we apply a similar comparison on the following linear media model of size $4\times 1$:
we take density $\rho$ that varies linearly in the range $[2,3]$, and Lam{\'e} parameters that varies linearly in the ranges $4\leq\lambda\leq 20$ and $1\leq\mu\leq 15$. The results, summarized in \ref{tab:2Dlinear}, are very similar to the results in Table \ref{tab:2Dhomogeneous}.

\begin{table} 
\centering
\begin{tabular}{c|ccc|ccc}
\toprule
  \mc{7}{c}{Iteration count for linear media}\\
\midrule
 &  \mc{3}{c|}{$\beta=1$, $G_s=10$} & \mc{3}{c}{$\beta=2/3$, $G_s=10 \, (8)$} \\
\midrule
 & \small 2-level & \small 3-level & \small 4-level & \small 2-level  & \small 3-level & \small 4-level \\
Grid size & {\small $\alpha = 0.1$} &  {\small $\alpha=0.4$}  & \small $\alpha=0.5$   & \small $\alpha=0.1$   & \small $\alpha=0.3 \,(0.4)$ & \small $\alpha=0.4 \,(0.5)$\\
\midrule
$512\times 128$ &   39 &  134  & 169  &  24 (31) &  77 (134) & 133  (181) \\
$1024\times 256$ &  59 &  265  & 333 &  34 (40) &  134 (219) &  260 (337) \\
$2048\times 512$ &  80 &  396  & 452  &  42 (52) &   153 (266) &  379 (479) \\
\bottomrule
 \end{tabular}
\caption{Number of preconditioning cycles needed for convergence with shifted Laplacian multigrid for the 2D elastic Helmholtz equation with the linear model in domain size $4\times1$. The frequency corresponds the number of grid-points per shear wavelength, $G_s$.}
\label{tab:2Dlinear}
\end{table}

\begin{figure}
  \centering
  \includegraphics[width=11cm]{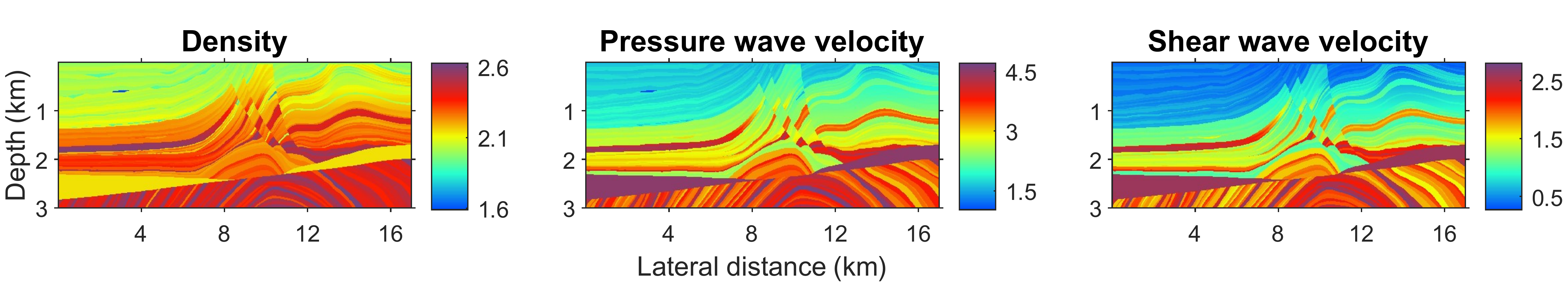}\\
  \caption{ The elastic Marmousi 2D model. Velocity units: $km/sec$, density units: $g/cm^3$.}\label{fig:Marmousi}
\end{figure}

Next, we give results of a similar comparison for the Marmousi-2 model \cite{martin2006marmousi2}. This geophysical 2D model is considered as case study for real 3D ground models. Since the model is very shallow --- only 3 $km$ depth, we add an extension of $0.5\, km$ on the bottom, to accommodate the absorbing boundary layer. We observe in Table \ref{tab:2Dmarmousi} that for two-grid cycles, the performance of our method with $\beta=2/3$ with $G_s=8$ is comparable $\beta=1$ with $G_s=10$, except for the 4-level method in the largest grid.

\begin{table} 
\centering
\begin{tabular}{c|ccc|ccc}
\toprule
  \mc{7}{c}{Iteration count for Marmousi-2 elastic media}\\
\midrule
 &  \mc{3}{c|}{$\beta=1$, $G_s=10$} & \mc{3}{c}{$\beta=2/3$, $G_s=10 \, (8)$} \\
\midrule
 & \small 2-level & \small 3-level & \small 4-level & \small 2-level  & \small 3-level & \small 4-level \\
Grid size & {\small $\alpha = 0.1$} &  {\small $\alpha=0.4$}  & \small $\alpha=0.5$   & \small $\alpha=0.1$   & \small $\alpha=0.3 \,(0.4)$ & \small $\alpha=0.4 \,(0.5)$\\
\midrule
$544\times 112$ & 31 & 107 & 147 &  28 (32) & 67 (106) & 100 (194) \\
$1088\times 224$ & 48 & 172 & 238  &  39 (51) & 102 (174) & 157 (338) \\
$2176\times 448$ & 70 & 243 & 374 & 53 (63) & 137 (242) & 238 (572) \\
\bottomrule
 \end{tabular}
\caption{Number of preconditioning cycles needed for convergence with shifted Laplacian multigrid for the 2D elastic Helmholtz equation with Marmousi-2 model. The frequency corresponds the number of grid-points per shear wavelength, $G_s$.}
\label{tab:2Dmarmousi}
\end{table}

\subsection{3D experiments} \label{subsec:3D}

\begin{figure}
  \centering
  \includegraphics[height=2cm]{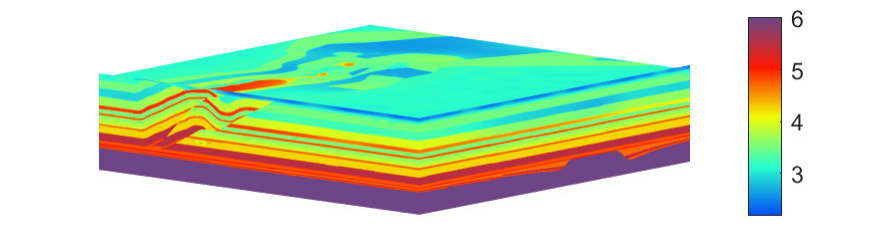}\\
  \caption{\footnotesize The SEG Overthrust pressure velocity model ($V_p$), in units of km/sec. The model corresponds to a domain of size 20$\times$20$\times$4.65 km.}\label{fig:Overthrust}
\end{figure}

In this subsection we give results for applying our system in 3D on the Overthrust velocity model described in Figure \ref{fig:Overthrust}. This is originally an acoustic model, and includes only pressure wave velocity, we define $V_s = 0.5V_p$ and $\rho = 0.25V_p + 1.5$. To implement the absorbing boundary layer, we add 16 grid points at the bottom of the domain.

Before presenting the results, we briefly describe the generalization of our discretization \eqref{eq:mixed-discretized} for the 3D case. The generalization is not straightforward, as there are various options how to spread a derivative in one direction, say, $x_1$, to the other two directions, $x_2,x_3$. The technique we chose here seems to replicate the behavior of the spread stencil in 2D, but further investigation is still needed via dedicated LFA and is beyond the scope of this paper. We define $\beta$-spread stencil for the first derivatives, for instance in the $x_1$ direction, as:
\begin{equation}\label{dxSkew3D}
(\partial_{x_1})^\beta_{h/2} = 
\beta
\begin{bmatrix}
-1 & * & 1
\end{bmatrix}
+(1-\beta)
\frac{1}{16}
\begin{bmatrix}
\begin{bmatrix}
-1 &   & 1 \\
-2 &   & 2 \\
-1 &   & 1 
\end{bmatrix}
&
\begin{bmatrix}
-2 &   & 2 \\
-4 & * & 4 \\
-2 &   & 2
\end{bmatrix}
&
\begin{bmatrix}
-1 &   & 1 \\
-2 &   & 2 \\
-1 &   & 1 
\end{bmatrix}
\end{bmatrix}.
\end{equation}
That is, we spread the $x_1$ derivative on both other directions $x_2$ and $x_3$. The derivatives in those directions are defined by stencils which are rotated versions of this stencil. Similarly to \eqref{eq:massHO}, the spread mass is defined by the stencil
\begin{equation} \label{eq:massHO3D}
M^\beta=
\omega^2 (1-\gamma \im) \left(
\beta\begin{bmatrix}
1 
\end{bmatrix}
+(1-\beta)\frac{1}{6}\begin{bmatrix}
\begin{bmatrix}
&  & \\
 & 1 & \\
&  & 
\end{bmatrix}
&
\begin{bmatrix}
& 1 & \\
1 & & 1 \\
& 1 & 
\end{bmatrix}
& 
\begin{bmatrix}
&  & \\
 & 1 & \\
&  & 
\end{bmatrix}
\end{bmatrix}
\right).
\end{equation}
The rest of the generalization is straightforward, where the main block in \eqref{eq:mixed-discretized} is a block diagonal matrix with three acoustic Helmholtz operators on its diagonal (rather than two).

Since this is a system of equations, the linear systems are huge even at rather high mesh sizes, requiring a lot of memory. Our implementation is applied in mixed percision, where the top level FGMRES method is applied using double precision, and the multigrid preconditioner is applied in single precision. Furthermore, the coarsest grid solution, even using 3 or 4 levels, is a challenging task. Here, we use the hybrid Kaczmarz iterative method described in subsection \ref{subsec:cycle}, where we apply the solver until a residual drop of 0.1. More specifically, we use FGMRES(5), where each preconditioning is applied using 10 parallel hybrid Kaczmarz steps using a damping of 0.8 and 8 cores. We limit the solution to be at most 250 Kaczmarz iterations. Better handling of this coarsest grid solution is part of our future research. Here we use K-cycles to accelerate the solution and maximize each inexact coarse solution.

Table \ref{tab:3Doverthrust} shows the 3D results. Our method, even down to 8 grid points per wavelength, performs significantly better than the standard discretization with 10 grid points per wavelength.

\begin{table} 
\centering
\begin{tabular}{c|cc|cc}
\toprule
  \mc{5}{c}{Iteration count for 3D Overthrust elastic media}\\
\midrule
 &  \mc{2}{c|}{$\beta=1$, $G_s=10$} & \mc{2}{c}{$\beta=2/3$, $G_s=10 \, (8)$} \\
\midrule
 &  \small 3-level & \small 4-level &  \small 3-level & \small 4-level \\
Grid size  &  {\small $\alpha=0.4$}  & \small $\alpha=0.5$   & \small $\alpha=0.3 \,(0.4)$ & \small $\alpha=0.4 \,(0.5)$\\
\midrule
$128\times128\times 56$ & 39 & 72 & 28 (30) &   31 (44) \\
$256\times 256\times 96$ & 77 &  175 & 49 (54) &   51 (70) \\
$384\times 384\times 136$ & 121 & $>$200 & 63 (67) &  73 (114) \\
\bottomrule
 \end{tabular}
\caption{Number of preconditioning cycles needed for convergence with shifted Laplacian multigrid for the 3D elastic Helmholtz equation with Overthrust model. The frequency corresponds the number of grid-points per shear wavelength, $G_s$.}
\label{tab:3Doverthrust}
\end{table}

\section{Conclusions} \label{sec:conclusion}
We introduced a novel discretization and a matrix-free multigrid method for the elastic Helmholtz equation in mixed formulation. We showed that our discretization, whose weights are tuned by two-grid LFA, is more suitable for shifted Laplacian geometric multigrid, yielding better performance compared to the standard finite-differences discretization for MAC staggered grid. We demonstrated that our discretization reduces numerical dispersion, especially on coarse grids. Our LFA results show that for 10 grid points per wavelength, the optimal weights in terms of multigrid convergence for the elastic Helmholtz equation in mixed formulation, coincides with the optimal weights in the acoustic version in terms of discretization error. We showed, numerically and theoretically, that our discretization allows the use of 8 grid points per wavelength, yielding at least the same performance as the standard stencil does for 10 grid points per wavelength. The stencil-based coarsening within the multigrid cycle makes our method easy to implement in a matrix-free form, hence suitable for parallel CPU and GPU computations.

\bibliographystyle{siamplain}
\bibliography{GeometricElasticHelmholtzSISC_ArxivV2.bbl}

\begin{thebibliography}{10}

\bibitem{aghamiry2022accurate}
{\sc H.~S. Aghamiry, A.~Gholami, L.~Combe, and S.~Operto}, {\em Accurate 3d
  frequency-domain seismic wave modeling with the wavelength-adaptive 27-point
  finite-difference stencil: A tool for full-waveform inversion}, Geophysics,
  87 (2022), pp.~R305--R324.

\bibitem{azulay2022multigrid}
{\sc Y.~Azulay and E.~Treister}, {\em Multigrid-augmented deep learning
  preconditioners for the helmholtz equation}, SIAM Journal on Scientific
  Computing,  (2022).

\bibitem{babuvska1995generalized}
{\sc I.~Babu{\v{s}}ka, F.~Ihlenburg, E.~T. Paik, and S.~A. Sauter}, {\em A
  generalized finite element method for solving the helmholtz equation in two
  dimensions with minimal pollution}, Computer methods in applied mechanics and
  engineering, 128 (1995), pp.~325--359.

\bibitem{baumann2018msss}
{\sc M.~Baumann, R.~Astudillo, Y.~Qiu, E.~Y. Ang, M.~Van~Gijzen, and R.-{\'E}.
  Plessix}, {\em An {MSSS}-preconditioned matrix equation approach for the
  time-harmonic elastic wave equation at multiple frequencies}, Computational
  Geosciences, 22 (2018), pp.~43--61.

\bibitem{bayliss1985accuracy}
{\sc A.~Bayliss, C.~I. Goldstein, and E.~Turkel}, {\em On accuracy conditions
  for the numerical computation of waves}, J. Comput. Phys., 59 (1985),
  pp.~396--404.

\bibitem{berenger1994perfectly}
{\sc J.-P. Berenger}, {\em A perfectly matched layer for the absorption of
  electromagnetic waves}, Journal of computational physics, 114 (1994),
  pp.~185--200.

\bibitem{bernth2011comparison}
{\sc H.~Bernth and C.~Chapman}, {\em A comparison of the dispersion relations
  for anisotropic elastodynamic finite-difference grids}, Geophysics, 76
  (2011), pp.~WA43--WA50.

\bibitem{Julia}
{\sc J.~Bezanson, A.~Edelman, S.~Karpinski, and V.~B. Shah}, {\em Julia: A
  fresh approach to numerical computing}, SIAM Review, 59 (2017), pp.~65--98,
  \url{https://doi.org/10.1137/141000671},
  \url{http://julialang.org/publications/julia-fresh-approach-BEKS.pdf},
  \url{https://arxiv.org/abs/http://dx.doi.org/10.1137/141000671}.

\bibitem{bonev2022hierarchical}
{\sc B.~Bonev and J.~S. Hesthaven}, {\em A hierarchical preconditioner for wave
  problems in quasilinear complexity}, SIAM Journal on Scientific Computing, 44
  (2022), pp.~A198--A229.

\bibitem{borisov20183d}
{\sc D.~Borisov, R.~Modrak, F.~Gao, and J.~Tromp}, {\em 3d elastic
  full-waveform inversion of surface waves in the presence of irregular
  topography using an envelope-based misfit function3d elastic fwi using
  envelopes}, Geophysics, 83 (2018), pp.~R1--R11.

\bibitem{brandt1977multi}
{\sc A.~Brandt}, {\em Multi-level adaptive solutions to boundary-value
  problems}, Mathematics of computation, 31 (1977), pp.~333--390.

\bibitem{briggs2000multigrid}
{\sc W.~L. Briggs, V.~E. Henson, and S.~F. McCormick}, {\em A multigrid
  tutorial}, SIAM, second~ed., 2000.

\bibitem{calandra2013improved}
{\sc H.~Calandra, S.~Gratton, X.~Pinel, and X.~Vasseur}, {\em An improved
  two-grid preconditioner for the solution of three-dimensional {Helmholtz}
  problems in heterogeneous media}, Numerical Linear Algebra with Applications,
  20 (2013), pp.~663--688.

\bibitem{cools2014new}
{\sc S.~Cools, B.~Reps, and W.~Vanroose}, {\em A new level-dependent coarse
  grid correction scheme for indefinite {Helmholtz} problems}, Numerical Linear
  Algebra with Applications, 21 (2014), pp.~513--533.

\bibitem{cools2013local}
{\sc S.~Cools and W.~Vanroose}, {\em Local {Fourier} analysis of the complex
  shifted {Laplacian} preconditioner for {Helmholtz} problems}, Numerical
  Linear Algebra with Applications, 20 (2013), pp.~575--597.

\bibitem{Dwarka2020}
{\sc V.~Dwarka and C.~Vuik}, {\em Scalable convergence using two-level
  deflation preconditioning for the helmholtz equation}, SIAM Journal on
  Scientific Computing, 42 (2020), pp.~A901--A928.

\bibitem{elec2001acoustic}
{\sc J.-C.~N. ELEC}, {\em Acoustic and electromagnetic equations: Integral
  representations for harmonic problems, vol. 44 of applied mathematical
  sciences}, 2001.

\bibitem{engquist1977absorbing}
{\sc B.~Engquist and A.~Majda}, {\em Absorbing boundary conditions for
  numerical simulation of waves}, Proceedings of the National Academy of
  Sciences, 74 (1977), pp.~1765--1766.

\bibitem{erlangga2006novel}
{\sc Y.~A. Erlangga, C.~W. Oosterlee, and C.~Vuik}, {\em A novel multigrid
  based preconditioner for heterogeneous {Helmholtz} problems}, SIAM J. Sci.
  Comput., 27 (2006), pp.~1471--1492.

\bibitem{gander2013domain}
{\sc M.~J. Gander and H.~Zhang}, {\em Domain decomposition methods for the
  {Helmholtz} equation: a numerical investigation}, in Domain Decomposition
  Methods in Science and Engineering XX, Springer, 2013, pp.~215--222.

\bibitem{gaspar2008distributive}
{\sc F.~Gaspar, J.~Gracia, F.~Lisbona, and C.~Oosterlee}, {\em Distributive
  smoothers in multigrid for problems with dominating grad--div operators},
  Numerical linear algebra with applications, 15 (2008), pp.~661--683.

\bibitem{gordon2013robust}
{\sc D.~Gordon and R.~Gordon}, {\em Robust and highly scalable parallel
  solution of the {Helmholtz} equation with large wave numbers}, Journal of
  Computational and Applied Mathematics, 237 (2013), pp.~182--196.

\bibitem{gosselin20143d}
{\sc B.~Gosselin-Cliche and B.~Giroux}, {\em {3D} frequency-domain
  finite-difference viscoelastic-wave modeling using weighted average 27-point
  operators with optimal coefficients}, Geophysics, 79 (2014), pp.~T169--T188.

\bibitem{greif2022block}
{\sc C.~Greif and Y.~He}, {\em Block preconditioners for the marker and cell
  discretization of the stokes-darcy equations}, arXiv preprint
  arXiv:2208.12357,  (2022).

\bibitem{gu201321}
{\sc B.~Gu, G.~Liang, and Z.~Li}, {\em A 21-point finite difference scheme for
  2d frequency-domain elastic wave modelling}, Exploration Geophysics, 44
  (2013), pp.~156--166.

\bibitem{he2022parameter}
{\sc Y.~He and Y.~Li}, {\em Parameter-robust braess-sarazin-type smoothers for
  linear elasticity problems}, arXiv preprint arXiv:2204.10462,  (2022).

\bibitem{heikkola2019parallel}
{\sc E.~Heikkola, K.~Ito, and J.~Toivanen}, {\em A parallel domain
  decomposition method for the helmholtz equation in layered media}, SIAM
  Journal on Scientific Computing, 41 (2019), pp.~C505--C521.

\bibitem{huismann2016fast}
{\sc I.~Huismann, J.~Stiller, and J.~Fr{\"o}hlich}, {\em Fast static
  condensation for the helmholtz equation in a spectral-element
  discretization}, in Parallel Processing and Applied Mathematics, Springer,
  2016, pp.~371--380.

\bibitem{jakobsen2020convergent}
{\sc M.~Jakobsen, R.-S. Wu, and X.~Huang}, {\em Convergent scattering series
  solution of the inhomogeneous helmholtz equation via renormalization group
  and homotopy continuation approaches}, Journal of Computational Physics, 409
  (2020), p.~109343.

\bibitem{kormann2017acceleration}
{\sc J.~Kormann, J.~E. Rodr{\'\i}guez, M.~Ferrer, A.~Farr{\'e}s,
  N.~Guti{\'e}rrez, J.~de~la Puente, M.~Hanzich, and J.~M. Cela}, {\em
  Acceleration strategies for elastic full waveform inversion workflows in {2D}
  and {3D}}, Computational Geosciences, 21 (2017), pp.~31--45.

\bibitem{levander1988fourth}
{\sc A.~R. Levander}, {\em Fourth-order finite-difference p-sv seismograms},
  Geophysics, 53 (1988), pp.~1425--1436.

\bibitem{li2022target}
{\sc Y.~Li and T.~Alkhalifah}, {\em Target-oriented high-resolution elastic
  full-waveform inversion with an elastic redatuming method}, Geophysics, 87
  (2022), pp.~R379--R389.

\bibitem{li2016Fourth}
{\sc Y.~Li, B.~Han, L.~M{\'e}tivier, and R.~Brossier}, {\em Optimal
  fourth-order staggered-grid finite-difference scheme for {3D}
  frequency-domain viscoelastic wave modeling}, J. Comput. Phys., 321 (2016),
  pp.~1055--1078.

\bibitem{li20152d}
{\sc Y.~Li, L.~M{\'e}tivier, R.~Brossier, B.~Han, and J.~Virieux}, {\em {2D}
  and {3D} frequency-domain elastic wave modeling in complex media with a
  parallel iterative solver}, Geophysics, 80 (2015), pp.~T101--T118.

\bibitem{maclachlan2011local}
{\sc S.~P. Maclachlan and C.~W. Oosterlee}, {\em A local fourier analysis
  framework for finite-element discretizations of systems of pdes}, Numer. Lin.
  Alg. Appl., 18 (2011), pp.~751--774.

\bibitem{martin2006marmousi2}
{\sc G.~S. Martin, R.~Wiley, and K.~J. Marfurt}, {\em Marmousi2: An elastic
  upgrade for marmousi}, The Leading Edge, 25 (2006), pp.~156--166.

\bibitem{mckee2008mac}
{\sc S.~McKee, M.~F. Tom{\'e}, V.~G. Ferreira, J.~A. Cuminato, A.~Castelo,
  F.~Sousa, and N.~Mangiavacchi}, {\em The mac method}, Computers \& Fluids, 37
  (2008), pp.~907--930.

\bibitem{notay2008recursive}
{\sc Y.~Notay and P.~S. Vassilevski}, {\em Recursive {Krylov}-based multigrid
  cycles}, Numerical Linear Algebra with Applications, 15 (2008), pp.~473--487.

\bibitem{operto20073d}
{\sc S.~Operto, J.~Virieux, P.~Amestoy, J.-Y. L’Excellent, L.~Giraud, and
  H.~B.~H. Ali}, {\em 3d finite-difference frequency-domain modeling of
  visco-acoustic wave propagation using a massively parallel direct solver: A
  feasibility study}, Geophysics, 72 (2007), pp.~SM195--SM211.

\bibitem{patera1984spectral}
{\sc A.~T. Patera}, {\em A spectral element method for fluid dynamics: laminar
  flow in a channel expansion}, Journal of computational Physics, 54 (1984),
  pp.~468--488.

\bibitem{pratt1999}
{\sc R.~Pratt}, {\em Seismic waveform inversion in the frequency domain, part
  1: Theory, and verification in a physical scale model}, Geophysics, 64
  (1999), pp.~888--901.

\bibitem{rizzuti2016multigrid}
{\sc G.~Rizzuti and W.~A. Mulder}, {\em Multigrid-based shifted-laplacian
  preconditioning for the time-harmonic elastic wave equation}, J. Comput.
  Phys., 317 (2016), pp.~47--65.

\bibitem{jInv17}
{\sc L.~Ruthotto, E.~Treister, and E.~Haber}, {\em {jInv} -- a flexible {Julia}
  package for {PDE} parameter estimation}, SIAM J. Sci. Comput., 39 (2017),
  pp.~S702–--S722.

\bibitem{saad1993flexible}
{\sc Y.~Saad}, {\em A flexible inner-outer preconditioned {GMRES} algorithm},
  SIAM J. Sci. Comput., 14 (1993), pp.~461--469.

\bibitem{singer1998high}
{\sc I.~Singer and E.~Turkel}, {\em High-order finite difference methods for
  the {Helmholtz} equation}, Computer Methods in Applied Mechanics and
  Engineering, 163 (1998), pp.~343--358.

\bibitem{vstekl1998accurate}
{\sc I.~{\v{S}}tekl and R.~G. Pratt}, {\em Accurate viscoelastic modeling by
  frequency-domain finite differences using rotated operators}, Geophysics, 63
  (1998), pp.~1779--1794.

\bibitem{stolk2016dispersion}
{\sc C.~C. Stolk}, {\em A dispersion minimizing scheme for the 3-d helmholtz
  equation based on ray theory}, Journal of computational Physics, 314 (2016),
  pp.~618--646.

\bibitem{stolk2014multigrid}
{\sc C.~C. Stolk, M.~Ahmed, and S.~K. Bhowmik}, {\em A multigrid method for the
  helmholtz equation with optimized coarse grid corrections}, SIAM Journal on
  Scientific Computing, 36 (2014), pp.~A2819--A2841.

\bibitem{sun2021deep}
{\sc H.~Sun and L.~Demanet}, {\em Deep learning for low-frequency extrapolation
  of multicomponent data in elastic fwi}, IEEE Transactions on Geoscience and
  Remote Sensing, 60 (2021), pp.~1--11.

\bibitem{sutmann2007compact}
{\sc G.~Sutmann}, {\em Compact finite difference schemes of sixth order for the
  helmholtz equation}, Journal of Computational and Applied Mathematics, 203
  (2007), pp.~15--31.

\bibitem{JointEikFWI17}
{\sc E.~Treister and E.~Haber}, {\em Full waveform inversion guided by travel
  time tomography}, SIAM J. Sci. Comput., 39 (2017), pp.~S587--–S609.

\bibitem{TREISTER2024112622}
{\sc E.~Treister and R.~Yovel}, {\em A hybrid shifted laplacian multigrid and
  domain decomposition preconditioner for the elastic helmholtz equations},
  Journal of Computational Physics, 497 (2024), p.~112622,
  \url{https://doi.org/https://doi.org/10.1016/j.jcp.2023.112622},
  \url{https://www.sciencedirect.com/science/article/pii/S0021999123007179}.

\bibitem{trottenberg2000multigrid}
{\sc U.~Trottenberg, C.~Oosterlee, and A.~Sch\"{u}ller}, {\em Multigrid},
  Academic Press, London and San Diego, 2001.

\bibitem{turkel2013compact}
{\sc E.~Turkel, D.~Gordon, R.~Gordon, and S.~Tsynkov}, {\em Compact {2D} and
  {3D} sixth order schemes for the {Helmholtz} equation with variable wave
  number}, J. Comput. Phys., 232 (2013), pp.~272--287.

\bibitem{umetani2009multigrid}
{\sc N.~Umetani, S.~P. MacLachlan, and C.~W. Oosterlee}, {\em A multigrid-based
  shifted laplacian preconditioner for a fourth-order {Helmholtz}
  discretization}, Numerical Linear Algebra with Applications, 16 (2009),
  pp.~603--626.

\bibitem{vanka1986blockFlow}
{\sc S.~Vanka}, {\em Block-implicit multigrid calculation of two-dimensional
  recirculating flows}, Computer Methods in Applied Mechanics and Engineering,
  59 (1986), pp.~29--48.

\bibitem{virieux2009overview}
{\sc J.~Virieux and S.~Operto}, {\em An overview of full-waveform inversion in
  exploration geophysics}, Geophysics, 74 (2009), pp.~WCC1--WCC26.

\bibitem{wang2020taylor}
{\sc B.~Wang, D.~Chen, B.~Zhang, W.~Zhang, M.~H. Cho, and W.~Cai}, {\em Taylor
  expansion based fast multipole method for 3-d helmholtz equations in layered
  media}, Journal of Computational Physics, 401 (2020), p.~109008.

\bibitem{wang2010acoustic}
{\sc S.~Wang, V.~Maarten, and J.~Xia}, {\em Acoustic inverse scattering via
  helmholtz operator factorization and optimization}, Journal of Computational
  Physics, 229 (2010), pp.~8445--8462.

\bibitem{wang20113d}
{\sc S.~Wang, V.~Maarten, and J.~Xia}, {\em On 3d modeling of seismic wave
  propagation via a structured parallel multifrontal direct helmholtz solver},
  Geophysical Prospecting, 59 (2011), pp.~857--873.

\bibitem{wobker2009numerical}
{\sc H.~Wobker and S.~Turek}, {\em Numerical studies of vanka-type smoothers in
  computational solid mechanics}, Advances in Applied Mathematics and
  Mechanics, 1 (2009), pp.~29--55.

\bibitem{zhu2010efficient}
{\sc Y.~Zhu, E.~Sifakis, J.~Teran, and A.~Brandt}, {\em An efficient multigrid
  method for the simulation of high-resolution elastic solids}, ACM
  Transactions on Graphics, 29 (2010), p.~16.

\bibitem{zubeldia2012energy}
{\sc M.~Zubeldia}, {\em Energy concentration and explicit sommerfeld radiation
  condition for the electromagnetic helmholtz equation}, Journal of Functional
  Analysis, 263 (2012), pp.~2832--2862.

\end{thebibliography}
\end{document}